%% file: main.tex
\newif\ifjpaa
\newif\ifarxiv
\jpaafalse\arxivtrue

\pdfoutput=1

\ifarxiv
\documentclass{article}
\else
\documentclass{elsarticle}
\fi

\usepackage{amsmath, amsthm, amssymb}
\usepackage{eucal}
\usepackage{geometry}
\usepackage[unicode,bookmarksnumbered=true,colorlinks=true,allcolors=blue,linktoc=all]{hyperref}
\usepackage[capitalise,nameinlink]{cleveref}
\usepackage{quiver}
\usepackage{tikz-cd}
\usepackage{stmaryrd}
\usepackage{listings}
\usepackage{etoolbox}
\usepackage{thm-restate}
\usepackage[nottoc,notlot,notlof]{tocbibind}
\usepackage{enumitem}
\usepackage{mathtools}
\ifarxiv
\usepackage{tocloft}
\fi

\crefname{subsection}{Subsection}{subsections}

\hypersetup{bookmarksdepth=2}
\newtheorem{theorem}{Theorem}

\newtheorem{thm}{Theorem}[section]
\newtheorem{lem}[thm]{Lemma}
\newtheorem{prop}[thm]{Proposition}
\newtheorem{cor}[thm]{Corollary}
\newtheorem{question}[thm]{Question}

\newtheorem{ssthm}{Theorem}[subsection]
\newtheorem{sslem}[ssthm]{Lemma}
\newtheorem{ssprop}[ssthm]{Proposition}

\theoremstyle{definition}
\newtheorem{defn}[thm]{Definition}

\newtheorem{ssdefn}[ssthm]{Definition}

\theoremstyle{remark}
\newtheorem{remark}[thm]{Remark}
\newtheorem{example}[thm]{Example}

\newtheorem{sswarn}[ssthm]{Warning}

\newcommand{\ncmd}{\newcommand}

\definecolor{DefColor}{rgb}{0.6,0.15,0.25}
\newcommand{\mdef}[1]{\textcolor{DefColor}{#1}}
\newcommand{\tdef}[1]{\mdef{\emph{#1}}}

\ncmd{\mbb}[1]{\mathbb{#1}}
\ncmd{\mrm}[1]{\mathrm{#1}}
\ncmd{\mcl}[1]{\mathcal{#1}}
\ncmd{\mfk}[1]{\mathfrak{#1}}
\ncmd{\mbf}[1]{\mathbf{#1}}

\usepackage{cjhebrew}
\DeclareFontFamily{U}{rcjhbltx}{}
\DeclareFontShape{U}{rcjhbltx}{m}{n}{<->s*[1.2]rcjhbltx}{}
\DeclareSymbolFont{hebrewletters}{U}{rcjhbltx}{m}{n}
\DeclareMathSymbol{\tsadi}{\mathord}{hebrewletters}{118}

\ncmd{\todo}[1]{\textbf{TODO #1}}
\ncmd{\reftodo}[1]{\textbf{REF #1}}

\ncmd{\id}[1][]{\mrm{id}\ifstrempty{#1}{}{_{#1}}}

\ncmd{\BB}{\mrm{B}}
\ncmd{\unit}{\mathbf{1}}
\ncmd{\CC}{\mcl{C}}
\ncmd{\DD}{\mcl{D}}
\ncmd{\EE}{\mcl{E}}
\ncmd{\bbE}{\mbb{E}}
\ncmd{\VV}{\mcl{V}}
\ncmd{\WW}{\mcl{W}}
\ncmd{\MM}{\mcl{M}}
\ncmd{\OO}{\mcl{O}}
\ncmd{\LM}{\mcl{LM}}
\ncmd{\Ass}{\mrm{Assoc}}
\ncmd{\II}{I}

\ncmd{\op}{\mrm{op}}
\ncmd{\rev}{\mrm{rev}}
\ncmd{\atomicpres}{\mrm{at}}
\ncmd{\atomic}[1][]{\ifstrempty{#1}{}{#1\text{-}}\atomicpres}
\ncmd{\pt}{\mrm{pt}}
\ncmd{\Hei}{\mathrm{Heine}}
\ncmd{\lax}{\mathrm{lax}}
\ncmd{\laxV}{{\lax\text{-}\VV}}
\ncmd{\cl}{\mathrm{cl}}

\ncmd{\LMod}{\mrm{LMod}}
\ncmd{\LModt}{\mbf{LMod}}
\ncmd{\LModtlax}{\LModt^\lax}
\ncmd{\wLMod}{\omega\LMod}
\ncmd{\wLModt}{\omega\LModt}
\ncmd{\LMODt}{\widehat{\mbf{LMod}}\vphantom{\LModt}}
\ncmd{\LMODtlax}{\LMODt^\lax}
\ncmd{\wLMODt}{\omega\LMODt}
\ncmd{\RMod}{\mrm{RMod}}
\ncmd{\Spaces}{\mcl{S}}
\ncmd{\SPACES}{\widehat{\Spaces}}
\ncmd{\Sp}{\mrm{Sp}}
\ncmd{\Fin}{\mrm{Fin}}
\ncmd{\Cat}{\mrm{Cat}}
\ncmd{\PCat}{\mrm{PCat}}
\ncmd{\CAT}{\widehat{\Cat}}
\ncmd{\CATcc}{\CAT\vphantom{\Cat}^{\mrm{cc}}}
\ncmd{\Catt}{\mbf{Cat}}
\ncmd{\CATt}{\widehat{\Catt}}
\ncmd{\CatV}{\Cat^\VV}
\ncmd{\CattV}{\Catt^\VV}
\ncmd{\CATV}{\CAT\vphantom{\Cat}^\VV}
\ncmd{\CATtV}{\CATt\vphantom{\Catt}^\VV}
\ncmd{\tCat}{\Cat_2}
\ncmd{\wCATV}{\omega\CATV}
\renewcommand{\Pr}{\mrm{Pr}}
\ncmd{\Prt}{\mbf{Pr}}
\ncmd{\LL}{\mrm{L}}
\ncmd{\RR}{\mrm{R}}
\ncmd{\all}{\mrm{all}}
\ncmd{\PrL}{\Pr^\LL}
\ncmd{\PrtL}{\Prt^\LL}
\ncmd{\PrR}{\Pr^\RR}
\ncmd{\PrLV}{\PrL_\VV}
\ncmd{\PrtLV}{\PrtL_\VV}
\ncmd{\intL}{\mrm{iL}}
\ncmd{\PriL}{\Pr^\intL}
\ncmd{\PrtiL}{\Prt^\intL}
\ncmd{\PriLV}{\PriL_\VV}
\ncmd{\PrtiLV}{\PrtiL_\VV}
\ncmd{\pwL}{{\LL\text{-}\mrm{pw}}}
\ncmd{\fwL}{{\LL\text{-}\mrm{fw}}}

\ncmd{\bfA}{\mbf{A}}
\ncmd{\bfX}{\mbf{X}}
\ncmd{\bfY}{\mbf{Y}}
\ncmd{\Cocart}{\mrm{Cocart}}
\ncmd{\Cocartt}{\mbf{Cocart}}
\ncmd{\CocarttO}{\Cocartt_\OO}
\ncmd{\CocarttOlax}{\CocarttO^\lax}
\ncmd{\Mon}{\mbf{Mon}}
\ncmd{\MonO}{\Mon_\OO}
\ncmd{\MonOlax}{\MonO^\lax}
\ncmd{\Funt}{\mbf{Fun}}
\ncmd{\Funtlax}{\Funt^\lax}

\ncmd{\yon}{\text{\usefont{U}{min}{m}{n}\symbol{'110}}}
\DeclareFontFamily{U}{min}{}
\DeclareFontShape{U}{min}{m}{n}{<-> dmjhira}{}
\ncmd{\yonV}{\yon{}^\VV}
\ncmd{\HeiyonV}{\yonV_\Hei}

\ncmd{\iso}{\xrightarrow{\smash{\raisebox{-0.4ex}{\ensuremath{\scriptstyle\sim}}}}}

\ncmd{\irchi}[2]{\raisebox{\depth}{$#1\chi$}}
\DeclareRobustCommand{\rchi}{{\mathpalette\irchi\relax}}

\DeclareMathOperator{\free}{free}
\DeclareMathOperator{\forget}{forget}
\DeclareMathOperator{\eval}{eval}
\DeclareMathOperator{\Fun}{Fun}
\DeclareMathOperator{\Alg}{Alg}
\DeclareMathOperator*{\colim}{colim}
\DeclareMathOperator{\coliminl}{colim}
\DeclareMathOperator{\PSh}{\mcl{P}}
\DeclareMathOperator{\PShV}{\PSh^\VV}
\DeclareMathOperator{\HeiPShV}{\PSh^\VV_\Hei}

\newcommand\noloc{%
  \nobreak
  \mspace{6mu plus 1mu}
  {:}
  \nonscript\mkern-\thinmuskip
  \mathpunct{}
  \mspace{2mu}
}

\ifjpaa
\begin{document}
	\begin{frontmatter}
\fi

\title{Naturality of the \texorpdfstring{$\infty$}{∞}-Categorical Enriched Yoneda Embedding}
\ifarxiv
\author{Shay Ben-Moshe\thanks{Einstein Institute of Mathematics, Hebrew University of Jerusalem.}}
\else
\author{Shay Ben-Moshe}
\ead{shay.benmoshe@mail.huji.ac.il}
\affiliation{organization={Einstein Institute of Mathematics, The Hebrew University of Jerusalem},
postcode={91904},
city={Jerusalem},
country={Israel}}
\fi
\date{}

\ifarxiv
\begin{document}
	\maketitle
\fi
	
	\begin{abstract}
		We make Hinich's $\infty$-categorical enriched Yoneda embedding natural.
		To do so, we exhibit it as the unit of a partial adjunction between the functor taking enriched presheaves and Heine's functor taking a tensored category to an enriched category.
		Furthermore, we study a finiteness condition of objects in a tensored category called being atomic, and show that the partial adjunction restricts to a (non-partial) adjunction between taking enriched presheaves and taking atomic objects.
	\end{abstract}

	\ifjpaa
	\begin{keyword}
		Enriched $\infty$-categories \sep Enriched Yoneda \sep Enriched presheaves \sep Atomic objects
		\MSC[2020] 18D20 \sep 18N70 \sep 18N60 \sep 18M05 \sep 18D25
	\end{keyword}
	\end{frontmatter}
	\fi
	
	\tableofcontents
	
	\input{intro.tex}
	\input{enrch.tex}
	\input{part-adj.tex}
	\input{2-cat.tex}
	\input{atomics.tex}
	\input{hei.tex}
	
	\bibliographystyle{alpha}
	\bibliography{refs}

\end{document}

%% file: intro.tex
\section{Introduction}

\subsection{Overview}

\subsubsection{Partial Adjunction}

The main goal of this paper is to make Hinich's enriched Yoneda embedding in the $\infty$-categorical setting from \cite{Hin} into a natural transformation.
We work in the framework of enriched $\infty$-categories as developed in \cite{Hin, Hin2, GH, Hei}.
These sources differ in notation, and we introduce our notation throughout the paper.
For brevity, henceforth we use the term (enriched) category to mean an (enriched) $\infty$-category, and a $2$-category to mean an $(\infty, 2)$-category.
Throughout this paper we fix a presentably monoidal category $\VV \in \Alg(\PrL)$.

For every $\VV$-enriched category $\CC_0 \in \CatV$ and a presentably $\VV$-tensored category $\DD \in \LMod_\VV(\PrL)$, Hinich defined an (unenriched) category of $\VV$-functors $\Fun_\VV(\CC_0, \DD) \in \PrL$.
Using this, he defined a presentably $\VV$-tensored category of enriched presheaves $\PShV(\CC_0) := \Fun_{\VV^\rev}(\CC_0^\op, \VV) \in \LMod_\VV(\PrL)$.
One of the main results of \cite{Hin} is the construction of a $\VV$-enriched Yoneda embedding
$$
\yonV\colon \CC_0 \to \PShV(\CC_0),
$$
namely, an object of $\Fun_\VV(\CC_0, \PShV(\CC_0)) \in \PrL$, satisfying the enriched Yoneda lemma.
However, this map is not shown to be natural in $\CC_0$, which is the central question of the present paper.

In a sequel paper \cite{Hin2}, Hinich shows that his enriched Yoneda embedding enjoys the following universal property:

\begin{thm}[{\cite[6.4.4]{Hin2}}]\label{hin-yon-intro}
    Let $\CC_0 \in \CatV$.
    Then, for every $\DD \in \LMod_\VV(\PrL)$, composition with the enriched Yoneda embedding induces an equivalence
    $$
    (\yonV)^*\colon
    \Fun^\LL_\VV(\PShV(\CC_0), \DD)
    \iso \Fun_\VV(\CC_0, \DD)
    .
    $$
\end{thm}

This almost exhibits $\yonV$ as a unit of an adjunction, but there are two problems:

\begin{enumerate}[ref=Problem (\arabic*)]
    \item\label{not-hom} $\Fun_\VV(\CC_0, \DD)$ is not the $\hom$ in any category (note that $\CC_0$ is $\VV$-enriched while $\DD$ is $\VV$-tensored).
    \item\label{size-issue} $\CC_0$ is small while $\DD$ is large.
\end{enumerate}

To solve \ref{not-hom}, we use the fact that presentably $\VV$-tensored categories produce $\VV$-enriched categories, as was first proven in \cite[Corollary 7.4.9]{GH}, and made functorial in the work of Heine \cite{Hei}.
Letting $\CATV$ denote the (huge) category of large $\VV$-enriched categories, a special case of Heine's result reads as follows:

\begin{thm}[{\cite[Corollary 6.13]{Hei}}]\label{chi-intro}
	There is a functor
	$$
	\rchi\colon \LMod_\VV(\PrL) \to \CATV
	$$
	witnessing the source as a (non-full non-wide) subcategory of the target.
\end{thm}

We shall also use \cite[Theorem 1.10]{Hei}, showing that for every $\CC_0 \in \CatV$ and $\DD \in \LMod_\VV(\PrL)$ there is an equivalence
$\Fun_\VV(\CC_0, \DD)^\simeq \cong \hom_\VV(\CC_0, \rchi(\DD))$,
where the right hand side denotes the space of morphisms in $\CATV$.
Combining these result, in \cref{psh-chi-pre-adj} we conclude that there is an isomorphism of spaces
$$
    \hom^\LL_\VV(\PShV(\CC_0), \DD)
    \iso \hom_\VV(\CC_0, \rchi(\DD)),
$$
where the left hand side denotes the underlying space of $\Fun^\LL_\VV(\PShV(\CC_0), \DD)$.

To deal with \ref{size-issue}, we recall that one can define adjoints partially, namely an adjoint and a (co)unit map defined only on a full subcategory (or, more generally, relative to a functor), see \cref{part-adj-def}.
Furthermore, we recall the folklore result that (partial) adjoints can be constructed point-wise, as we show in \cref{part-adj}.

With these results in mind, we deduce our first main result:

\begin{theorem}[{\cref{pshv-chi-adj}}]\label{pshv-chi-adj-intro}
	The functor $\rchi\colon \LMod_\VV(\PrL) \to \CATV$ has a partial left adjoint $\PShV\colon \CatV \to \LMod_\VV(\PrL)$ with partial unit agreeing with the enriched Yoneda embedding.
\end{theorem}

In particular, for every $f\colon \CC_0 \to \DD_0$ in $\CatV$, we get an induced $\VV$-linear left adjoint $f_!\colon \PShV(\CC_0) \to \PShV(\DD_0)$, and an isomorphism $f_! \yonV \cong \yonV f$, as explained in \cref{yon-nat}.

\subsubsection{Atomics--Presheaves Adjunction}

In a somewhat different direction, we study atomic objects, a finiteness condition on objects of presentably $\VV$-tensored categories.
We show that this condition is closely related to the enriched Yoneda embedding, and, in particular, gives another solution for \ref{size-issue}, leading to a (non-partial) adjunction.

Let $\CC \in \LMod_\VV(\PrL)$, and let $X \in \CC$.
The functor $- \otimes X\colon \VV \to \CC$ is a $\VV$-linear left adjoint functor.
We say that $X$ is \emph{atomic} if the right adjoint $\hom^\VV(X, -)$ is itself a left adjoint and the canonical lax $\VV$-linear structure on it is strong (see \cref{at-def}).
These objects form a small full $\VV$-subcategory $\CC^{\atomic} \subset \rchi(\CC)$.
We show that the construction $\CC \mapsto \CC^{\atomic}$ is functorial in \emph{internally left adjoint functors}, that is, $\VV$-linear left adjoint functors $L\colon \CC \to \DD \in \LMod_\VV(\PrL)$ whose right adjoint is itself a left adjoint and the canonical lax $\VV$-linear structure on it is strong.
Namely, we obtain a functor
$$
(-)^{\atomic}\colon \LMod_\VV(\PrL)^\intL \to \CatV.
$$

As a key example, in \cref{yon-at} we show that for any $\CC_0 \in \CatV$ and $X \in \CC_0$, the image $\yonV(X) \in \PShV(\CC_0)$ is atomic.
To see this, recall that by the enriched Yoneda lemma $\hom^\VV(\yonV(X), -)$ is given by evaluation at $X$, which preserves (co)limits and is $\VV$-linear, as we show in \cref{colim-levelwise}.
Therefore, we get a factorization of the enriched Yoneda embedding through the atomics $\yonV\colon \CC_0 \to \PShV(\CC_0)^{\atomic}$.
Furthermore, we show that the functoriality of $\PShV$ obtained in \cref{pshv-chi-adj-intro} sends $\VV$-functors to internally left adjoints, namely, restricts to a functor
$$
\PShV\colon \CatV \to \LMod_\VV(\PrL)^\intL.
$$
Our second main result is that the partial adjunction of \cref{pshv-chi-adj-intro} restricts accordingly:

\begin{theorem}[{\cref{yon-at-adj}}]\label{pshv-at-adj-intro}
	There is an adjunction
	$$
	\PShV\colon \CatV \rightleftarrows \LMod_\VV(\PrL)^\intL\noloc (-)^{\atomic}
	$$
	with unit agreeing with the enriched Yoneda embedding.
\end{theorem}

In particular, for every $f\colon \CC_0 \to \DD_0$ in $\CatV$, the $\VV$-linear left adjoint $f_!\colon \PShV(\CC_0) \to \PShV(\DD_0)$, is internally left adjoint, and thus admits a right adjoint $f^\circledast\colon \PShV(\DD_0) \to \PShV(\CC_0)$ which is itself a $\VV$-linear left adjoint.
Using this, in \cref{upperstar-pointwise} we show that $f^\circledast(G)(X) \cong G(f(X))$.
We warn the reader that this does not imply that $f^\circledast$ is the functor given by pre-composition, see \cref{nat-precomp} for further discussion.

\subsubsection{$2$-Categorical Structures}

In our discussion of atomic objects and their relation to the enriched Yoneda embedding, we use certain $2$-categorical aspects of the theory of enriched categories and tensored categories.
Notably, we have considered the condition for a $\VV$-linear left adjoint functor to be internally left adjoint, which we show has a particularly simple $2$-categorical interpretation.

In \cref{lmod-lax} we recall that the category of $\VV$-tensored categories and lax $\VV$-linear functors enhances to a $2$-category $\LModtlax_\VV$.
We show that a $1$-morphism in that $2$-category is a left adjoint, if and only if it is strong $\VV$-linear and left adjoint (on the underlying category).
This result is a special case of a more general result we deduce from Lurie's work on relative adjunctions and \cite{laxoplax}, concerning the $2$-category $\MonOlax$ of $\OO$-monoidal categories and lax $\OO$-monoidal functors for some operad $\OO$, which may be of independent interest.

\begin{prop}[{\cref{lax-O-adj}}]
	The left adjoints in $\MonOlax$ are the lax $\OO$-monoidal functors that are strong and fiber-wise left adjoint.
\end{prop}

The category $\LMod_\VV(\PrL)$ also enhances into a $2$-category, which by the above result is a full $2$-subcategory $\LModt_\VV(\PrtL) \subset (\LMODtlax_\VV)^\LL$ of the left adjoints in the large version of the $2$-category $\LModtlax_\VV$.
Furthermore, using the above result again, we see that internally left adjoints are precisely the left adjoint morphisms in the $2$-category $\LModt_\VV(\PrtL)$.

As explained in \cite[Construction 8.2]{Hei-arxiv}, by \cite[Proposition 5.7.16]{GH}, the category of $\VV$-enriched categories also enhances to a $2$-category $\CattV$.
Furthermore, by \cite[Theorem 8.3]{Hei-arxiv}, the functor $\rchi$ of \cref{chi-intro} enhances to a $2$-functor.
Using these results we also deduce the following result, which may be of independent interest, and is closely related to \cite[Theorem 1.2]{Hei}.
Here a $2$-functor is called \emph{$2$-fully faithful} if it induces an isomorphism on $\hom$ categories, i.e.\ if it exhibits the source as a full $2$-subcategory of the target.

\begin{theorem}[{\cref{chi-2}, \cref{chi-2-l}}]
	The functor $\rchi$ enhances to a $2$-fully faithful $2$-functor
	$$
	\rchi\colon \LModt_\VV(\PrtL) \to (\CATtV)^\LL.
	$$
    Taking left adjoints again we get the $2$-fully faithful $2$-functor
	$$
	\rchi\colon \LModt_\VV(\PrtL)^\intL \to (\CATtV)^{\LL\LL}.
	$$
\end{theorem}

In particular, for $f\colon \CC_0 \to \DD_0$ in $\CatV$, since $f_!\colon \PShV(\CC_0) \to \PShV(\DD_0)$ is internally left adjoint, we get a corresponding double adjunction $\rchi(f_!) \dashv \rchi(f^\circledast) \dashv \rchi(f_\circledast)$ in $\CATtV$.

\begin{remark}
    The $2$-categorical enhancement of $\rchi$ mentioned above is not incorporated in the published version of Heine's paper \cite{Hei}.
    However, it does appear in a newer (and in an older) version of the paper available on arXiv \cite{Hei-arxiv}, which we cite for this result.
\end{remark}

\subsection{Relation to Previous Work}

This paper shares many of the ideas on atomic objects developed in our previous paper with Tomer Schlank, most notably \cite[Theorem D]{BMS}.
To avoid the usage of enriched categories, in our previous paper we work over a \emph{mode} $\MM$, that is, an idempotent algebra in $\PrL$.
The key feature of modes is that $\LMod_\MM(\PrL)$ is a \emph{full} subcategory of $\PrL$, namely, a left adjoint functor is automatically (and uniquely) $\MM$-linear.
This allowed us, for example, to simplify the definition of $X$ being atomic to having $\hom^\MM(X, -)\colon \CC \to \MM$ commute with colimits (and thus automatically $\MM$-linear).
Particularly, in the previous paper we have ignored the $\MM$-enriched structure on the category of atomic objects.
In addition, we worked with unenriched presheaves and the unenriched Yoneda embedding, further tensored with $\MM$.
In particular, \cite[Theorem D]{BMS} is a weaker statement than the main results of the present paper in these respects.

On the other hand, in the unenriched context, the two different functorialities of presheaves and naturalities of the Yoneda embedding were shown to coincide, contrary to the enriched case of the present paper (see the discussion in \cref{nat-precomp}).
In addition, our goal in the previous paper was different.
The unenriched presheaves functor is endowed with a symmetric monoidal structure, and through the adjunction, this makes the atomics functor lax symmetric monoidal, and the Yoneda a symmetric monoidal natural transformation.
In the present paper we do not deal with the multiplicative structure (though see the discussion in \cref{mul}).

\subsection{Further Questions}\label{subsec-quest}

We now list several further questions left open, which we expect have a positive answer, but we do not know how to approach.

\subsubsection{Naturality for Pre-composition}\label{nat-precomp}

The first question is closely related to \cite{Ram} and \cite[Corollary F]{HHLN}, which deal with the case $\VV = \Spaces$.

Recall that using the universal property of $\PShV(\CC_0)$, we deduced \cref{pshv-chi-adj-intro}, assembling $\PShV$ into a functor, and the enriched Yoneda embedding into a natural transformation.
In addition, we saw that $f_!$ has a right adjoint $f^\circledast$, and we showed that $f^\circledast(G)(X) \cong G(f(X))$.

The construction of $\PShV(\CC_0)$ as enriched functors from $\CC_0^\op$ to $\VV$ shows that it admits another functoriality in $\CC_0$.
More specifically, \cite[6.1.4]{Hin} implies that it assembles into a functor $\PShV\colon \CatV \to (\PrR)^\op$, sending $f$ to $f^*\colon \PShV(\DD_0) \to \PShV(\CC_0)$ given by pre-composition, admitting a left adjoint $f_?$, commonly spelled ``$f$ lower what".
Note that this functoriality does not a priori give the $\VV$-tensored structure.
Also note that $f^\circledast(G)(X) \cong G(f(X)) \cong f^*(G)(X)$, however it is not clear that this holds naturally in $X$, $G$ or $\CC_0$.

Furthermore, Hinich's construction of the enriched Yoneda embedding, as described for example in \cite[8.4.1]{Hin2}, seems to interact with the $f^*$, and thus $f_?$, functoriality, but we don't know how to extract naturality from his results.

\begin{question}
    Can the $f_?$ functoriality be extended to $\VV$-modules?
    Can the enriched Yoneda embedding be made natural for the $f_?$ functoriality?
    Do these agree with the $f_!$ functoriality and naturality of the enriched Yoneda embedding?
\end{question}

\subsubsection{Heine's and Hinich's Enriched Yoneda Embeddings}\label{hei-hin}

In \cite{Hei}, Heine also defines an enriched Yoneda embedding $\HeiyonV\colon \CC_0 \to \rchi(\PShV(\CC_0))$ by different means.
As we explain in \cref{sec-hei}, the main results of this paper hold for this version as well, producing an adjunction
$$
\HeiPShV\colon \CatV \rightleftarrows \LMod_\VV(\PrL)^\intL\noloc (-)^{\atomic}
$$
with unit agreeing with Heine's enriched Yoneda embedding.
The uniqueness of adjoints implies that there is a natural isomorphism $\psi\colon \PShV \iso \HeiPShV$, together with an isomorphism $\psi \yonV \cong \HeiyonV$.
On the other hand, by construction, point-wise we have $\HeiPShV(\CC_0) = \PShV(\CC_0)$, so we get an automorphism $\psi_{\CC_0}\colon \PShV(\CC_0) \to \HeiPShV(\CC_0) = \PShV(\CC_0)$.
It is not clear to us that this automorphism is the identity.
Showing this is equivalent to showing that these two versions of the enriched Yoneda embedding coincide.

\begin{question}\label{hin-hei-same}
    Do the enriched Yoneda embeddings constructed by Hinich and Heine coincide?
    Namely, is the above automorphism the identity?
\end{question}

\subsubsection{$\CAT$-Enrichment}

Recall that $\rchi$ enhances to a $2$-functor.
Furthermore, the universal property of the enriched Yoneda embedding of \cref{hin-yon-intro} is originally stated for functor categories, though one of the sides is not constructed as the $\CAT$-enriched $\hom$ in a category.
This leads to the following question:

\begin{question}
    Can the (partial) adjunctions be made $\CAT$-enriched?
\end{question}

\subsubsection{Multiplicative Structure}\label{mul}

For simplicity, assume that $\VV$ is presentably symmetric monoidal.
\cite[Theorem 7.19]{Hei} enhances $\rchi$ into a symmetric monoidal functor.
Via the adjunction, it may be possible to endow $\PShV$ with an oplax symmetric monoidal structure, either by making the subfunctor $(-)^{\atomic}$ lax symmetric monoidal and using the main result of \cite{laxoplax}, or by proving a version of it for partial adjunctions and directly applying to $\PShV$.
Analogously to the case of $\VV = \Spaces$, one would expect the resulting oplax symmetric monoidal structure on $\PShV$ to be strong, and thus make the adjunction symmetric monoidal.
Assuming this, for any operad $\OO$, we get that if $\CC_0$ is $\OO$-monoidal, then $\PShV(\CC_0)$ and the enriched Yoneda embedding are endowed with an $\OO$-monoidal structure.

Separately, in \cite[8.3 and 8.4]{Hin2}, Hinich studies the multiplicative structure of enriched presheaves and the enriched Yoneda embedding.
In particular, under the above assumptions, he endows $\PShV(\CC_0)$ and the enriched Yoneda embedding with an $\OO$-monoidal structure.
He further proves an $\OO$-monoidal version of the universal property of \cref{hin-yon-intro}.

\begin{question}
    Can the adjunction be made symmetric monoidal?
    Does the induced $\OO$-monoidal structure on $\PShV(\CC_0)$ and the enriched Yoneda embedding agree with those constructed by Hinich?
\end{question}

\subsection{Acknowledgments}

I would like to thank Lior Yanovski for numerous useful conversations about atomic objects, enriched categories and weighted colimits.
I would also like thank Hadrian Heine for generously answering several questions regarding his work, for updating \cite{Hei} to include some results used in the present paper and for writing \cite{Hei2} upon my question on the topic.
I also thank Shai Keidar, Shaul Ragimov, Maxime Ramzi and Noam Zimhoni for comments on earlier drafts of this paper.
Particularly, I would like to thank Maxime for pointing several subtle missing components in earlier drafts, and suggesting corrections for some of them.
Finally, I thank the anonymous referee for their careful reading and helpful comments.

%% file: enrch.tex
\section{Generalities on Enriched Categories}\label{sec-enrch}

In this section we review some generalities on enriched categories, their relationship to tensored categories, and the enriched Yoneda embedding.
We shall not delve into the details of the constructions, as most of them will not play a role in the present paper, but rather only the formal properties of the resulting objects.

\subsection{Enriched Categories and Tensored Categories}

\begin{defn}
	We denote the (large) \tdef{category of $\VV$-enriched categories}, defined in \cite[7.1.2]{Hin}, by $\mdef{\CatV}$.
	For $\CC_0, \DD_0 \in \CatV$, we denote the space of \tdef{$\VV$-functors} between them by $\mdef{\hom_\VV(\CC_0, \DD_0)} \in \Spaces$.
	Similarly, we let $\mdef{\CATV}$ be the (huge) category of large $\VV$-enriched categories.
\end{defn}

$\VV$-enriched categories are closely related to categories tensored over $\VV$, as was first proven in \cite[Corollary 7.4.9]{GH}, and made functorial in the work of Heine.
Indeed, Heine constructs the category $\mdef{\wLMod^\lax_\VV}$ (denoted $\wLMod_\VV$ there) of \tdef{weakly $\VV$-tensored categories and lax $\VV$-linear functors}, and a full subcategory $\mdef{\wLMod^{\cl,\lax}_\VV}$ thereof on the \emph{closed} weakly $\VV$-tensored categories.

\begin{thm}[{\cite[Theorem 1.5]{Hei}}]\label{chi-w}
	There is an equivalence
	$$
	\mdef{\rchi}\colon \wLMod^{\cl,\lax}_\VV \iso \CatV.
	$$
\end{thm}

Considering the large version of this equivalence, one can restrict the source to the subcategory with objects presentable categories with $\VV$-action commuting with colimits (which are automatically closed) and morphisms the (strong) $\VV$-linear left adjoint functors.

\begin{defn}\label{prlv}
	We define the category of \tdef{presentably $\VV$-tensored categories} to be
	$\mdef{\PrLV} := \LMod_\VV(\PrL)$.
	For $\CC, \DD \in \PrLV$, we denote the space of \tdef{$\VV$-linear left adjoint functors} between them by $\mdef{\hom^\LL_\VV(\CC, \DD)} \in \SPACES$, which is the space of objects of the category $\mdef{\Fun^\LL_\VV(\CC, \DD)} \in \CAT$.
\end{defn}

\begin{cor}[{\cite[Corollary 6.13]{Hei}}]\label{chi}
	There is a functor
	$$
	\mdef{\rchi}\colon \PrLV \to \CATV.
	$$
	witnessing the source as a (non-full non-wide) subcategory of the target.
\end{cor}

\begin{remark}
	In \cite[Theorem 1.2]{Hei} Heine uses the notation $\PrLV$ for the image of the inclusion above.
\end{remark}

For $\CC \in \PrLV$, this constructs $\rchi(\CC) \in \CATV$, both of which have the same underlying category and thus space of objects.

\subsection{Enriched Yoneda Lemma and Weighted Colimits}

In \cite{Hin, Hin2} Hinich constructs enriched presheaves and the enriched Yoneda embedding, which we recall in this subsection.

We begin with Hinich's model for the category of $\VV$-functors from a $\VV$-enriched category to a $\VV$-tensored category.
In Hinich's model, a $\VV$-enriched category $\CC_0 \in \CatV$ is an algebra in some operad constructed from the space of objects $\CC_0^{\simeq}$ (see \cite[3.1.1]{Hin}).
For a presentably $\VV$-tensored category $\DD \in \PrLV$, Hinich endows the (unenriched) category of functors $\Fun(\CC_0^\simeq, \DD)$ with a left module structure over this operad (see \cite[6.1.1]{Hin}).
In particular, one can consider left $\CC_0$-modules, arriving at the following definition.

\begin{defn}[{\cite[6.1.3]{Hin}}]\label{fun-model}
	Let $\CC_0 \in \CatV$ and $\DD \in \PrLV$, then the category of \tdef{$\VV$-functors} from $\CC_0$ to $\DD$ is defined to be
	$$
	\mdef{\Fun_\VV(\CC_0, \DD)} := \LMod_{\CC_0}(\Fun(\CC_0^\simeq, \DD))
	\in \PrL.
	$$
\end{defn}

We note that $\Fun_\VV(\CC_0, \DD)$ is indeed presentable.
To see this, note that $\Fun(\CC_0^\simeq, \DD)$ is an (unenriched) presheaf category, and thus is presentable by \cite[Proposition 5.5.3.6]{HTT}.
Then, the module category $\Fun_\VV(\CC_0, \DD)$ is also presentable by \cite[Corollary 4.2.3.7]{HA}.

\begin{defn}[{\cite[6.2.2]{Hin}}]\label{pshv-def}
	Let $\CC_0 \in \CatV$.
	The $\VV$-tensored category of \tdef{$\VV$-enriched presheaves} is defined to be
	$$
	\mdef{\PShV(\CC_0)} := \Fun_{\VV^\rev}(\CC_0^\op, \VV) \in \PrLV.
	$$
\end{defn}

\begin{remark}
	$\VV$ is a $\VV$-$\VV$-bimodule in $\PrL$.
	One of the $\VV$-module structures is used to define the presentable category of $\VV^\rev$-enriched functors, and the other is used to endow the resulting category with a $\VV$-module structure.
\end{remark}

\begin{remark}
	In \cite{Hei} Heine gives an a priori different, although very similar, definition of $\VV$-functors and $\VV$-enriched presheaves.
	Nevertheless, \cite{Hei2} shows that his definitions agree with Hinich's definition, and henceforth we use them interchangeably.
\end{remark}

We now record the main results about Hinich's enriched Yoneda embedding.
We remark that in \cite[6.2.7]{Hin}, Hinich proves his version of the enriched Yoneda lemma, very closely related to \ref{yon-lem} of \cref{hin-yon}.
However, the form described below is somewhat different, and relies on a definition of $\hom^\VV$ and the evaluation at $X$ that were not introduced thus far.
We postpone these definitions to \cref{hom-v-def} and \cref{eval-x}, and the proof to \cref{yon-lem-restate}, but include the statement here for completeness of the exposition.
Note that \ref{yon-lem} (as well as \ref{wcol-lin}) will not be used before \cref{sec-at}.

\begin{thm}[{\cite{Hin, Hin2}}]\label{hin-yon}
	Let $\CC_0 \in \CatV$.
	Then, there is a \tdef{$\VV$-enriched Yoneda embedding} $\VV$-functor
	$$
	\mdef{\yonV}\colon \CC_0 \to \PShV(\CC_0),
	$$
	that is, an object of $\Fun_\VV(\CC_0, \PShV(\CC_0))$.
	For every $\DD \in \PrLV$, there is a \tdef{weighted colimit} functor
	$$
	\mdef{\coliminl_{\CC_0}^{(-)}(-)}\colon \PShV(\CC_0) \times \Fun_\VV(\CC_0, \DD) \to \DD
	.
	$$
	These satisfy the following properties:
	\begin{enumerate}[ref=(\arabic*)]
		\item\label{yon-lem} For every $X \in \CC_0$, the functor $\hom^\VV(\yonV(X), -)\colon \PShV(\CC_0) \to \VV$ agrees with evaluation at $X$ as $\VV$-linear functors.
		\item\label{wcol-lin} $\coliminl_{\CC_0}^{(-)}(-)$ commutes with colimits in both arguments separately, and commutes with the $\VV$-action in the first argument.
		\item\label{psh-univ} There is an equivalence
			$$
			(\yonV)^*\colon
			\Fun^\LL_\VV(\PShV(\CC_0), \DD)
			\rightleftarrows \Fun_\VV(\CC_0, \DD)
			\noloc \coliminl_{\CC_0}^{(-)}(-)
			.
			$$
	\end{enumerate}
\end{thm}

\begin{proof}
	$\yonV$ is constructed in \cite[6.2.4]{Hin}.
	\ref{yon-lem} is deferred to \cref{yon-lem-restate}, though see \cite[6.2.7]{Hin}.
	$\coliminl_{\CC_0}^{(-)}(-)$ is constructed in \cite[6.2.3]{Hin2}, where \ref{wcol-lin} is explained.
	\ref{psh-univ} is \cite[6.4.4]{Hin2}.
\end{proof}

Heine shows that $\Fun_\VV(\CC_0, \DD)$ is closely related to the $\hom$ in $\VV$-enriched categories via $\rchi$ (see also \cite[6.3.6]{Hin}).
We recall a special case of Heine's result as follows:

\begin{prop}[{\cite[Theorem 1.10]{Hei}}]\label{fun-model-enrch}
	Let $\CC_0 \in \CatV$ and $\DD \in \PrLV$, then there is an equivalence
	$$
	\Fun_\VV(\CC_0, \DD)^\simeq \cong \hom_\VV(\CC_0, \rchi(\DD))
	$$
	natural in $\CC_0$ and $\DD$.
	Namely, there is a natural isomorphism between the functors
	$$
	(\CatV)^\op \times \PrLV \xrightarrow{i \times \rchi} (\CATV)^\op \times \CATV \xrightarrow{\hom_\VV} \SPACES
	$$
	and
	$$
	(\CatV)^\op \times \PrLV \xrightarrow{\Fun_\VV(-, -)} \Cat \xrightarrow{(-)^\simeq} \Spaces \subset \SPACES
	$$
\end{prop}

We use this to transform the enriched Yoneda embedding to a morphism in $\CATV$, and deduce a universal property similar to \ref{psh-univ} of \cref{hin-yon} where both the source and the target are the $\hom$ in some category.

\begin{defn}\label{yon-chi}
	We denote by the same notation the enriched Yoneda embedding $\VV$-functor $\yonV\colon \CC_0 \to \rchi(\PShV(\CC_0))$, corresponding to the enriched Yoneda embedding under the equivalence
	$\Fun_\VV(\CC_0, \PShV(\CC_0))^\simeq \cong \hom_\VV(\CC_0, \rchi(\PShV(\CC_0)))$.
\end{defn}

\begin{prop}\label{psh-chi-pre-adj}
	Let $\CC_0 \in \CatV$ and $\DD \in \PrLV$, then the composition
	$$
	\hom^\LL_\VV(\PShV(\CC_0), \DD)
	\xrightarrow{\rchi} \hom_\VV(\rchi(\PShV(\CC_0)), \rchi(\DD))
	\xrightarrow{(\yonV)^*} \hom_\VV(\CC_0, \rchi(\DD))
	$$
	is an equivalence.
\end{prop}

\begin{proof}
	Naturality in the $\PrLV$ coordinate of \cref{fun-model-enrch} shows that for $\CC_0 \in \CatV$ and $\DD, \EE \in \PrLV$ we get a commutative square:
	$$\begin{tikzcd}
		\hom^\LL_\VV(\EE, \DD) \arrow{r}{\rchi} \arrow{d}{} & \hom_\VV(\rchi(\EE), \rchi(\DD)) \arrow{d}{} \\
		\hom(\Fun_\VV(\CC_0, \EE)^\simeq, \Fun_\VV(\CC_0, \DD)^\simeq) \arrow{r}{\sim} & \hom(\hom_\VV(\CC_0, \rchi(\EE)), \hom_\VV(\CC_0, \rchi(\DD)))
	\end{tikzcd}$$
	Using the exponential adjunction we get a commutative square:
	$$\begin{tikzcd}
		\hom^\LL_\VV(\EE, \DD) \times \Fun_\VV(\CC_0, \EE)^\simeq \arrow{r}{\rchi} \arrow{d}{\circ}
			& \hom_\VV(\rchi(\EE), \rchi(\DD)) \times \hom_\VV(\CC_0, \rchi(\EE)) \arrow{d}{\circ} \\
		\Fun_\VV(\CC_0, \DD)^\simeq \arrow{r}{\sim} & \hom_\VV(\CC_0, \rchi(\DD))
	\end{tikzcd}$$
	Taking $\EE = \PShV(\CC_0)$, and picking the point $\yonV \in \Fun_\VV(\CC_0, \PShV(\CC_0))^\simeq$, we get a commutative square:
	$$\begin{tikzcd}
		\hom^\LL_\VV(\PShV(\CC_0), \DD) \arrow{r}{\rchi} \arrow{d}{(\yonV)^*} & \hom_\VV(\rchi(\PShV(\CC_0)), \rchi(\DD)) \arrow{d}{(\yonV)^*} \\
		\Fun_\VV(\CC_0, \DD)^\simeq \arrow{r}{\sim} & \hom_\VV(\CC_0, \rchi(\DD))
	\end{tikzcd}$$
	The left morphism is an equivalence by applying $(-)^\simeq$ to \ref{psh-univ} of \cref{hin-yon}.
	Thus, the left-bottom composition is an equivalence, and by the commutativity of the diagram, so is the upper-right composition, concluding the proof.
\end{proof}

%% file: part-adj.tex
\section{Partial Adjunction}\label{sec-part-adj}

In this section we prove \cref{pshv-chi-adj-intro}, namely the naturality of the enriched Yoneda embedding.
To do so, we first define (unenriched) partial adjunctions.
Then, we show the folklore result that (partial) adjoints can be constructed point-wise, from which the naturality of the enriched Yoneda embedding follows immediately.

\begin{defn}\label{part-adj-def}
	Let $L\colon \CC_0 \to \DD_0$ be a functor, and $i\colon \EE_0 \to \DD_0$ another functor.
	The data of a \tdef{partial right adjoint} (of $L$ relative to $i$) is a functor $R\colon \EE_0 \to \CC_0$ and a natural transformation $\varepsilon\colon LR \Rightarrow i$ called the \tdef{partial counit}, such that for every $X \in \CC_0$ and $Y \in \EE_0$ the composition
	$$
	\hom_{\CC_0}(X, RY)
	\xrightarrow{L} \hom_{\DD_0}(LX, LRY)
	\xrightarrow{\varepsilon_Y \circ -} \hom_{\DD_0}(LX, iY)
	$$
	is an isomorphism.
	A partial left adjoint with a partial unit is defined dually, by taking $(-)^\op$.
\end{defn}

\begin{remark}
	What we call a partial adjunction is typically called a relative adjunction.
	However, in \cref{sec-2-cat} we use the distinct concept of relative adjunctions developed by Lurie.
	To avoid confusion, we call what is typically called a relative adjunction a partial adjunction.
\end{remark}

\begin{remark}
	In this paper we only use the notion of a partial adjunction when $i\colon \EE_0 \to \DD_0$ is an inclusion of a full subcategory.
	In this case, one can think of $R$ as an adjoint defined only partially on $\DD_0$, justifying the name.
\end{remark}

\begin{remark}
	Unlike in a standard adjunction, a partial adjunction is not a symmetric concept.
	Particularly, note that there is no partial unit (in fact, $L$ and $R$ can not be composed in the other direction).
\end{remark}

\begin{prop}\label{part-adj}
	Let $L\colon \CC_0 \to \DD_0$ be a functor, and $i\colon \EE_0 \to \DD_0$ another functor.
	Assume that for every $Y \in \EE_0$ we are given an object $RY \in \CC_0$ and a morphism $\varepsilon_Y\colon LRY \to iY$, such that for every $X \in \CC_0$ the composition
	$$
	\hom_{\CC_0}(X, RY)
	\xrightarrow{L} \hom_{\DD_0}(LX, LRY)
	\xrightarrow{\varepsilon_Y \circ -} \hom_{\DD_0}(LX, iY)
	$$
	is an isomorphism.
	Then, $R$ and $\varepsilon$ assemble into the data of a partial right adjoint.
\end{prop}

\begin{proof}
	Consider the functor $\tilde{R}\colon \DD_0 \to \PSh(\CC_0)$ given by the composition
	$$
	\DD_0 \xrightarrow{\yon_{\DD_0}} \PSh(\DD_0) \xrightarrow{L^*} \PSh(\CC_0),
	$$
	which sends $Z \in \DD_0$ to $\hom_{\DD_0}(L(-), Z)\colon \CC_0^\op \to \Spaces$.

	Let $Y \in \EE_0$.
	Consider the natural transformation
	\begin{equation}\label{L-eps}
		\hom_{\CC_0}(-, RY)
		\xrightarrow{L} \hom_{\DD_0}(L(-), LRY)
		\xrightarrow{\varepsilon_Y \circ -} \hom_{\DD_0}(L(-), iY)
		= \tilde{R}(iY).
	\end{equation}
	By assumption, it is an isomorphism at every $X \in \CC_0$, and thus it is a natural isomorphism.
	Namely, the presheaf $\tilde{R}(iY) \in \PSh(\CC_0)$ is representable by $RY \in \CC_0$.
	In other words, the composition
	$\tilde{R} i\colon \EE_0 \to \PSh(\CC_0)$
	lands in the essential image of the (unenriched) Yoneda embedding $\yon_{\CC_0}\colon \CC_0 \to \PSh(\CC_0)$.
	As the Yoneda embedding is fully faithful, we get an induced functor $R\colon \EE_0 \to \CC_0$, together with a natural isomorphism $\yon_{\CC_0} R \cong \tilde{R} i = L^* \yon_{\DD_0} i$ of functors $\EE_0 \to \PSh(\CC_0)$.

	By construction, $R$ agrees with $RY$ at every $Y \in \EE_0$.
	Furthermore, the natural isomorphism at $Y$ is given by \eqref{L-eps}, constructed via $\varepsilon_Y$.
	We now extract the partial counit $\varepsilon\colon LR \Rightarrow i$.
	Consider the following composition:
	\begin{align*}
		\yon_{\DD_0} LR
		&\overset{\text{(1)}}{\cong} (LR)_! \yon_{\EE_0}\\
		&\overset{\text{(2)}}{\Rightarrow} (LR)_! R^* \yon_{\CC_0} R\\
		&\overset{\text{(3)}}{\cong} (LR)_! R^* L^* \yon_{\DD_0} i\\
		&\cong (LR)_! (LR)^* \yon_{\DD_0} i\\
		&\overset{\text{(4)}}{\Rightarrow} \yon_{\DD_0} i.
	\end{align*}
	(1) is naturality of the (unenriched) Yoneda embedding.
	(2) is the fact that $R$ is a functor, giving $\hom_{\EE_0}(-, -) \Rightarrow \hom_{\CC_0}(R(-), R(-))$, which by exponential adjunction is the same as $\yon_{\EE_0} \Rightarrow R^* \yon_{\CC_0} R$.
	(3) is the isomorphism $\yon_{\CC_0} R \cong L^* \yon_{\DD_0} i$.
	(4) is the counit of the adjunction $(LR)_! \dashv (LR)^*$.
	Thus, we have constructed a natural transformation $\yon_{\DD_0} LR \Rightarrow \yon_{\DD_0} i$.
	Since $\yon_{\DD_0}$ is fully faithful, this induces a natural transformation $\varepsilon\colon LR \Rightarrow i$.
	Recall that the isomorphism used at step (3) is given at $Y \in \EE_0$ by \eqref{L-eps}, which shows that $\varepsilon$ indeed agrees with $\varepsilon_Y$.
\end{proof}

\begin{thm}\label{pshv-chi-adj}
	The functor $\rchi\colon \PrLV \to \CATV$ has a partial left adjoint $\PShV\colon \CatV \to \PrLV$
	with partial unit $\yonV\colon \id|_{\CatV} \to \rchi \PShV$,
	agreeing with enriched presheaves and enriched Yoneda embedding.
\end{thm}

\begin{proof}
	This follows immediately from \cref{psh-chi-pre-adj} and (the dual statement to) \cref{part-adj}.
\end{proof}

\begin{defn}\label{pshv-of-f}
	For $f\colon \CC_0 \to \DD_0$ a morphism in $\CatV$, we denote the induced morphism in $\PrLV$ by $\mdef{f_!}\colon \PShV(\CC_0) \to \PShV(\DD_0)$.
\end{defn}

An instance of the naturality of the enriched Yoneda embedding is the following:

\begin{prop}\label{yon-nat}
	Let $f\colon \CC_0 \to \DD_0$.
	Then, $\yonV_{\DD_0} f \cong \rchi(f_!) \yonV_{\CC_0}$ in $\hom_\VV(\CC_0, \rchi(\PShV(\DD_0)))$.
	Similarly, $\yonV_{\DD_0} f \cong f_! \yonV_{\CC_0}$ in $\Fun_\VV(\CC_0, \PShV(\DD_0))$ of \cref{fun-model}.
\end{prop}

\begin{proof}
	The first part is an immediate application of the adjunction of \cref{pshv-chi-adj}.
	For the second part, recall from \cref{fun-model-enrch} that there is a natural isomorphism $\hom_\VV(-, \rchi(-)) \cong \Fun_\VV(-, -)^\simeq$.
	Recall that in \cref{yon-chi} we have defined $\yonV_{\CC_0} \in \hom_\VV(\CC_0, \rchi(\PShV(\CC_0)))$ to be the map corresponding to $\yonV_{\CC_0}\in \Fun_\VV(\CC_0, \PShV(\CC_0))$ via this isomorphism, and similarly for $\DD_0$.
	By applying the naturality of the isomorphism in the target to the morphism $f_!$, we get that $\rchi(f_!) \yonV_{\CC_0} \in \hom_\VV(\CC_0, \rchi(\PShV(\DD_0)))$ corresponds to $f_! \yonV_{\CC_0} \in \Fun_\VV(\CC_0, \PShV(\DD_0))$.
	Similarly, by applying the naturality of the isomorphism in the source to the morphism $f$, we get that $\yonV_{\DD_0} f \in \hom_\VV(\CC_0, \rchi(\PShV(\DD_0)))$ corresponds to $\yonV_{\DD_0} f \in \Fun_\VV(\CC_0, \PShV(\DD_0))$.
	Thus, under the natural isomorphism, the isomorphism $\yonV_{\DD_0} f \cong \rchi(f_!) \yonV_{\CC_0}$ corresponds to an isomorphism $\yonV_{\DD_0} f \cong f_! \yonV_{\CC_0}$.
\end{proof}

%% file: 2-cat.tex
\section{\texorpdfstring{$2$}{2}-Categorical Structures}\label{sec-2-cat}

In this section we study certain $2$-categorical aspects of enriched categories and tensored categories.
The first main result of this section is \cref{lax-O-adj}, showing that a lax $\OO$-monoidal functor between $\OO$-monoidal categories is a left adjoint in the $2$-category $\MonOlax$, if and only if it is strong $\OO$-monoidal and fiber-wise left adjoint.
This is proven by applying Lurie's \cite[\S7.3.2 Relative Adjunctions]{HA} to the $2$-categories constructed in \cite{laxoplax}.
From this, we deduce \cref{lax-lin-adj}, the corresponding result for lax $\VV$-linear functors between $\VV$-tensored categories.
Then, after recalling that in \cite[Theorem 8.3]{Hei-arxiv} Heine shows that $\rchi$ is a $2$-functor, we deduce \cref{chi-2}, our second main result, which says that $\rchi$ enhances to a $2$-fully faithful $2$-functor
$\rchi\colon \PrtLV \to (\CATtV)^\LL$.
Finally, in \cref{yon-lem-restate} we finish the proof of \ref{yon-lem} of \cref{hin-yon}.

\setcounter{subsection}{-1}
\subsection{Left Adjoints in \texorpdfstring{$2$}{2}-Categories}\label{subsec-left-adj}

In the published version of this paper, the proof of \cref{lax-lin-adj} contains a slight oversight in showing that $(-)^\LL$ preserve limits.
This subsection was added to make the argument precise (note that the numberings in the rest of the paper were not changed).

\begin{ssdefn}
	For a $2$-category $\bfX$, we denote by $\mdef{\bfX^\LL} \subset \bfX$ the $1$-full $2$-subcategory on all objects, $1$-morphisms that are left adjoints, and all higher morphisms.
	Noting that $2$-functors send left adjoints to left adjoints, by \cite[Proposition A.1]{Ram} we get a subfunctor of the identity functor
	\[
		(-)^\LL\colon \tCat \to \tCat.
	\]
\end{ssdefn}

\begin{ssdefn}
	We define the \tdef{$2$-cell} to be the $2$-category
	\[\begin{tikzcd}
		{\mdef{C_2}} & 0 && 1
		\arrow["{:=}"{description}, draw=none, from=1-1, to=1-2]
		\arrow[""{name=0, anchor=center, inner sep=0}, curve={height=12pt}, from=1-2, to=1-4]
		\arrow[""{name=1, anchor=center, inner sep=0}, curve={height=-12pt}, from=1-2, to=1-4]
		\arrow[shorten <=3pt, shorten >=3pt, Rightarrow, from=1, to=0]
	\end{tikzcd}\]
	We denote by $\mdef{\partial C_2} \subset C_2$ the underlying category of $C_2$, namely two objects connected by two parallel $1$-morphisms.
\end{ssdefn}

\begin{sslem}\label{cells-jc}
	The functors
	\[
		\hom(\pt, -)\colon \tCat \to \Spaces,
		\quad \hom([1], -)\colon \tCat \to \Spaces,
		\quad \hom(C_2, -)\colon \tCat \to \Spaces
	\]
	are jointly conservative and limit preserving.
\end{sslem}

\begin{proof}
	These functors are corepresentable, hence limit preserving, and it remains to show that they are jointly conservative.
	Considering $\pt, [1], C_2$ as strict $2$-categories, they generate Joyal's category $\Theta_2$ under pushouts.
	By \cite{Camp}, it follows that they generate $\Theta_2$ under pushouts also as (weak) $2$-categories.
	Since $\Theta_2$ is dense in $\tCat$, it follows that the $2$-categories $\pt, [1], C_2$ generate $\tCat$ under colimits, and the result follows from \cite[Corollary 2.5]{MonTow}.
\end{proof}

\begin{ssprop}\label{L-lims}
	The functor $(-)^\LL$ preserves limits.
\end{ssprop}

\begin{proof}
	By \cref{cells-jc}, it suffices to show that the compoistions
	\[
		\hom(\pt, (-)^\LL)\colon \tCat \to \Spaces,
		\quad \hom([1], (-)^\LL)\colon \tCat \to \Spaces,
		\quad \hom(C_2, (-)^\LL)\colon \tCat \to \Spaces
	\]
	preserve limits.
	
	For $\pt$, note that the canonical map $\bfX^\LL \to \bfX$ induces an isomorphism on underlying spaces.
	In other words, $\hom(\pt, (-)^\LL) \simeq \hom(\pt, -)$, and the latter preserves limits since it is corepresentable.

	For $[1]$, note that $\hom([1], (-)^\LL)$ is the space of adjunctions.
	By \cite{Adj}, this functor is corepresented by the walking adjunction, hence it preserves limits.
	
	For $C_2$, recall that, by definition, $\bfX^\LL$ is a $1$-full $2$-subcategory of $\bfX$.
	Therefore, the following is a pullback square:
	\[\begin{tikzcd}
		{\hom(C_2,(-)^\LL)} & {\hom(C_2,-)} \\
		{\hom(\partial C_2,(-)^\LL)} & {\hom(\partial C_2,-)}
		\arrow[from=1-1, to=1-2]
		\arrow[from=1-1, to=2-1]
		\arrow[from=1-2, to=2-2]
		\arrow[from=2-1, to=2-2]
	\end{tikzcd}\]
	The functors at the top-right and bottom-right are corepresentable, hence preserve limits.
	Note that $\partial C_2$ is equivalent to $[1] \coprod_{[1]^\simeq} [1]$ as strict $2$-categories, thus, by \cite{Camp}, the same is true as (weak) $2$-categories.
	Moreover, $[1]^\simeq$ is equivalent to $\pt \coprod \pt$.
	Therefore, the cases of $\pt$ and $[1]$ imply that the functor at the bottom-left preserves limits.
	Since the top-left functor is the pullback of the other three, it follows that it preserves limits as well, concluding the proof.
\end{proof}

\begin{sswarn}
	We have not defined $(-)^\LL$ as a $2$-functor, but we expect that it does lift to a $2$-functor.
	However, such a $2$-functor will \emph{not} preserve $2$-limits (i.e.\ lax limits).
	To see this, consider the case of a constant diagram $\bfX\colon \II \to \tCat$ where $\II$ is a category.
	As constant $2$-limits are $2$-functor $2$-categories, in this case, the assembly map is
	$$
	\Funt(\II, \bfX)^\LL \to \Funt(\II, \bfX^\LL)
	$$
	which is an inclusion but not generally an equivalence, by \cite[Corollary 5.2.12]{fun-adj}.
\end{sswarn}

\subsection{Adjunctions in \texorpdfstring{$\OO$}{O}-Monoidal Categories}

\begin{defn}[{\cite[Definition 3.1.7]{laxoplax}}]
	Let $B$ be a a category.
	Denote by $\mdef{\Cocartt_B^\lax}$ the $2$-category of cocartesian fibrations over $B$ and functors over $B$.
\end{defn}

We shall not recall the precise definition of this $2$-category (for which we refer the reader to \cite[Definition 3.1.7]{laxoplax}), rather, let us recall its objects, $1$-morphisms and $2$-morphisms.

\begin{prop}
	$\Cocartt_B^\lax$ has
	\begin{itemize}
		\item objects: cocartesian fibrations $q\colon \CC \to B$,
		\item $1$-morphisms: functors over $B$, namely, a $1$-morphism from $q\colon \CC \to B$ to $p\colon \DD \to B$ is a functor $F\colon \CC \to \DD$ and a natural isomorphism $p F \cong q$,
		\item $2$-morphisms: natural transformations over $B$, namely, a $2$-morphism from $F$ to $G$ is a natural transformation $\alpha\colon F \Rightarrow G$ and an identification of the natural transformation $q \cong p F \xRightarrow{p \alpha} p G \cong q$ with $\id_q$, i.e.\ exhibiting the following square as commutative:
		\[\begin{tikzcd}
			pF & pG \\
			q & q
			\arrow["p\alpha", Rightarrow, from=1-1, to=1-2]
			\arrow["\wr", Rightarrow, from=1-2, to=2-2]
			\arrow["\wr"', Rightarrow, from=1-1, to=2-1]
			\arrow[Rightarrow, no head, from=2-1, to=2-2]
		\end{tikzcd}\]
	\end{itemize}
\end{prop}

\begin{proof}
	By construction, the underlying category of the $2$-category $\Cocartt_B^\lax$ is the full subcategory of $\Cat_{/B}$ on the cocartesian fibrations, which explains the objects and the $1$-morphisms.
	For the $2$-morphisms, we recall that as explained in \cite[Observation 3.1.9]{laxoplax}, a $2$-morphism is commutative diagram of the form:
	\[\begin{tikzcd}
		{\CC \times [1]} && \DD \\
		& B
		\arrow["q"', from=1-1, to=2-2]
		\arrow["p", from=1-3, to=2-2]
		\arrow["\bar\alpha", from=1-1, to=1-3]
	\end{tikzcd}\]
	Under the exponential adjunction, $\bar\alpha\colon \CC \times [1] \to \DD$ corresponds to the natural transformation $\alpha\colon F \Rightarrow G$, and the isomorphism $p \bar\alpha \cong q$ corresponds to the identification $p \alpha \cong \id_q$.
\end{proof}

\begin{lem}\label{2-inv}
	A $2$-morphism in $\Cocartt_B^\lax$ given by $\alpha\colon F \Rightarrow G$ together with an identification $p \alpha \cong \id_q$ is invertible if and only if $\alpha$ is invertible.
\end{lem}

\begin{proof}
	Clearly, if the $2$-morphism is invertible then in particular $\alpha$ is.
	For the other direction, assume that $\alpha$ has an inverse $\alpha^{-1}\colon G \Rightarrow F$.
	We shall enhance it to an inverse $2$-morphism in $\Cocartt_B^\lax$.
	Consider the following diagram:
	\[\begin{tikzcd}
		pF & pG & pF \\
		q & q & q
		\arrow["\wr"', Rightarrow, from=1-1, to=2-1]
		\arrow["\wr"', Rightarrow, from=1-3, to=2-3]
		\arrow["\wr"', Rightarrow, from=1-2, to=2-2]
		\arrow["p\alpha", Rightarrow, from=1-1, to=1-2]
		\arrow["{p\alpha^{-1}}", Rightarrow, from=1-2, to=1-3]
		\arrow[Rightarrow, no head, from=2-1, to=2-2]
		\arrow[Rightarrow, no head, from=2-2, to=2-3]
	\end{tikzcd}\]
	The identification $p \alpha \cong \id_q$ is precisely an (invertible) $3$-morphism in $\Catt$ making the left square commute.
	The outer square commutes as $p\alpha^{-1} p\alpha \cong \id_{pF} \cong \id_q$.
	Thus, by composing the outer square with the inverse of the left square, we get commutativity data for the right square.
	This makes $\alpha^{-1}$ into a $2$-morphism in $\Cocartt_B^\lax$, which by construction is the required inverse.
\end{proof}

Our next goal is to understand left adjoints in the $2$-category $\Cocartt_B^\lax$, achieved in \cref{cocart-l}.
To that end, we shall employ Lurie's \cite[\S7.3.2 Relative Adjunctions]{HA}.
The definitions and results in that section are phrased for cartesian fibrations and right adjoints, which by taking $(-)^\op$ correspond to our case of cocartesian fibrations and left adjoints, as we present here.

First, we recall that in a $2$-category, we say that a $1$-morphism $L\colon X \to Y$ is a \emph{left adjoint} if there exists a $1$-morphism $R\colon Y \to X$ and two $2$-morphisms $u\colon \id_X \Rightarrow RL$ and $c\colon LR \Rightarrow \id_Y$ satisfying the zigzag identities, namely
$$
(R \xRightarrow{uR} RLR \xRightarrow{Rc} R) \cong \id_R,
\quad
(L \xRightarrow{Lu} LRL \xRightarrow{cL} L) \cong \id_L.
$$
We also recall that to check that $L$ is a left adjoint, we can check a weaker condition:

\begin{lem}[{\cite[2.1.11]{RV}}]\label{min-adj}
	In any $2$-category, a $1$-morphism $L\colon X \to Y$ is a left adjoint if and only if there exist a $1$-morphism $R\colon Y \to X$ and $2$-morphisms $u\colon \id_X \Rightarrow RL$ and $c\colon LR \Rightarrow \id_Y$ such that the zigzag morphisms
	$$
	R \xRightarrow{uR} RLR \xRightarrow{Rc} R,
	\quad
	L \xRightarrow{Lu} LRL \xRightarrow{cL} L
	$$
	are \emph{invertible} (in which case there exists a possibly different $2$-morphism $\tilde{c}\colon LR \Rightarrow \id_Y$ for which the zigzag identities hold).
\end{lem}

\begin{proof}
	The proof of \cite[2.1.11]{RV} works in an arbitrary $2$-category, although stated in the context of functors between ($\infty$-)categories.
\end{proof}

\begin{defn}[{\cite[Definition 7.3.2.2]{HA}}]\label{rel-adj}
	Let
	\[\begin{tikzcd}
		\CC && \DD \\
		& B
		\arrow["q"', from=1-1, to=2-2]
		\arrow["p", from=1-3, to=2-2]
		\arrow["L", from=1-1, to=1-3]
	\end{tikzcd}\]
	be a commutative diagram of categories.
	We say that $L$ \tdef{admits a relative right adjoint} if $L$ has a (non-relative) right adjoint $R\colon \DD \to \CC$, and the counit map $c\colon LR \Rightarrow \id_{\DD}$ satisfies the condition that $pc\colon pLR \Rightarrow p$ is an equivalence.
\end{defn}

\begin{remark}
	The definition of \cite[Definition 7.3.2.2]{HA} is via the two equivalent conditions of \cite[Proposition 7.3.2.1]{HA}.
	First, we note that these are phrased for admitting a relative left adjoint, but the theory is symmetric via taking $(-)^\op$.
	Second, Lurie assumes that $p$ and $q$ are categorical fibrations, but every functor is equivalent to a categorical fibration, so this assumption can be dropped.
	Finally, and most importantly, we note that our definition clearly implies condition (1) and is implied by condition (2) in \cite[Proposition 7.3.2.1]{HA}, and thus is also equivalent to them.
\end{remark}

\begin{prop}\label{adj-rel-adj}
	A $1$-morphism $L\colon \CC \to \DD$ over $B$ in $\Cocartt_B^\lax$ is a left adjoint if and only if $L$ admits a relative right adjoint.
\end{prop}

\begin{proof}
	If $L$ is a left adjoint in $\Cocartt_B^\lax$ then it clearly admits a relative right adjoint.
	For the other direction, assume that $L$ admits a relative right adjoint, and we shall show that it is a left adjoint in $\Cocartt_B^\lax$.
	Let $R\colon \DD \to \CC$, $u\colon \id_\CC \Rightarrow RL$ and $c\colon LR \Rightarrow \id_\DD$ be the (non-relative) adjunction data, such that $pc$ is an equivalence.
	Note that the equivalence $pL \cong q$ and $pc$ together provide an equivalence $qR \cong pLR \cong p$, making $R$ into a $1$-morphism in $\Cocartt_B^\lax$.
	Moreover, the fact that $pc$ is an equivalence and the canonical commutativity of the square
	\[\begin{tikzcd}
		pLR & p \\
		p & p
		\arrow["{ }", Rightarrow, no head, from=2-1, to=2-2]
		\arrow[Rightarrow, no head, from=1-2, to=2-2]
		\arrow["pc", Rightarrow, from=1-1, to=1-2]
		\arrow["pc"', Rightarrow, from=1-1, to=2-1]
	\end{tikzcd}\]
	makes $c$ into a $2$-morphism from $LR$ to $\id_\DD$ in $\Cocartt_B^\lax$.
	In the rest of the proof we make $u$ into a $2$-morphism as well, and show that $L$ is a left adjoint using \cref{min-adj}.

	We begin by making $u$ into a $2$-morphism.
	Consider the following diagram:
	\[\begin{tikzcd}
		pL && pLRL && pL \\
		q && q && q
		\arrow["pLu", Rightarrow, from=1-1, to=1-3]
		\arrow["pcL", Rightarrow, from=1-3, to=1-5]
		\arrow["\wr"', Rightarrow, from=1-3, to=2-3]
		\arrow["\wr"', Rightarrow, from=1-1, to=2-1]
		\arrow["\wr"', Rightarrow, from=1-5, to=2-5]
		\arrow[Rightarrow, no head, from=2-1, to=2-3]
		\arrow[Rightarrow, no head, from=2-3, to=2-5]
	\end{tikzcd}\]
	The upper composition is $p(-)$ applied the zigzag morphism, which is equivalent to $\id_{pL}$, thus making the outer square commute.
	The right square canonically commutes.
	Thus, by composing the outer square with the inverse of the right square, we get commutativity data for the left square as well.
	In other words, we obtained an identification of $pLu$ with $\id_q$.
	The identification $pL \cong q$ thus shows that we got an identification of $qu$ with $\id_q$, making $u$ into a $2$-morphism.

	By \cref{min-adj}, to see that $L$ is a left adjoint in $\Cocartt_B^\lax$ it suffices to check that the zigzag morphisms corresponding to these $2$-morphisms are invertible.
	By \cref{2-inv}, a $2$-morphism in $\Cocartt_B^\lax$ is invertible if and only if the underlying natural transformations are invertible.
	Indeed, the underlying natural transformations are simply
	$$
	R \xRightarrow{uR} RLR \xRightarrow{Rc} R,
	\quad
	L \xRightarrow{Lu} LRL \xRightarrow{cL} L
	$$
	which are invertible (in fact, equivalent to the identities) by the assumption that $u$ and $c$ are a unit and a counit for a (non-relative) adjunction $L \dashv R$.
\end{proof}

The following is a recast of \cite[Proposition 7.3.2.6]{HA} in an appropriate $2$-category.

\begin{prop}\label{cocart-l}
	The left adjoints in $\Cocartt_B^\lax$ are those functors over $B$ that preserve cocartesian morphisms and are fiber-wise left adjoint
	$$
	(\Cocartt_B^\lax)^\LL = \Cocartt_B^\fwL.
	$$
\end{prop}

\begin{proof}
	By \cref{adj-rel-adj}, a $1$-morphism $L$ between two cocartesian fibrations is a left adjoint if and only if it admits a relative right adjoint.
	Thus by (the $(-)^\op$ of) \cite[Proposition 7.3.2.6]{HA}, $L$ is a left adjoint if and only if it preserves locally cocartesian morphisms and is fiber-wise left adjoint.
	By \cite[Proposition 2.4.2.8]{HTT}, locally cocartesian morphisms and cocartesian morphisms coincide, concluding the proof.
\end{proof}

\begin{remark}
	There is a possible alternative route to \cref{cocart-l}.
	Recall that \cite[Theorem E]{laxoplax} shows that $\Cocartt_B^\lax \cong \Funtlax(B, \Catt)$, where the right hand side denotes the $2$-category of functors and lax natural transformations.
	Additionally, \cite[Theorem 4.6]{laxnat} shows that for any $2$-categories $\bfX, \bfY$, the left adjoints in $\Funtlax(\bfX, \bfY)$ are those lax natural transformations which are strong and point-wise left adjoint.
	Combined, this would imply \cref{cocart-l}.
	However, the definition of $\Funtlax$ in \cite{laxoplax} and \cite{laxnat} rely on two different models of the Gray tensor product, which to the best of our knowledge were not shown to be equivalent so far, and thus this does not constitute a complete proof.
\end{remark}

In the rest of the subsection we follow the notations of \cite[3.4]{laxoplax} for operads.
Namely, an operad $\OO$ is a functor $\OO \to \Fin_*$ (satisfying certain properties), whereas in \cite{HA}, it is typically denoted by $\OO^\otimes$.

\begin{defn}[{\cite[Definition 3.4.1]{laxoplax}}]
	Let $\OO$ be an operad.
	Let $\mdef{\MonOlax} \subset \CocarttOlax$ be the $1$-full $2$-subcategory on the $\OO$-monoidal categories and lax $\OO$-monoidal functors (i.e.\ functors over $\OO$ preserving cocartesian morphisms lying over inert morphisms).
\end{defn}

The following is a recast of \cite[Corollary 7.3.2.7]{HA} in an appropriate $2$-category.

\begin{prop}\label{lax-O-adj}
	The left adjoints in $\MonOlax$ are those lax $\OO$-monoidal functors that are strong and fiber-wise left adjoint
	$$
	(\MonOlax)^\LL = \MonO^\fwL.
	$$
\end{prop}

\begin{proof}
	For the first direction, note that a left adjoint in $\MonOlax$ is in particular a left adjoint in $\CocarttOlax$, thus it is strong $\OO$-monoidal (preserves all cocartesian morphisms) and fiber-wise left adjoint.

	For the second direction, let $L\colon \CC \to \DD$ over $\OO$ in $\MonO^\fwL$.
	As $\MonO^\fwL \subset \CocarttO^\fwL = (\CocarttOlax)^\LL$, there is a right adjoint $R\colon \DD \to \CC$ over $\OO$ in $\CocarttOlax$ with some unit and counit.
	Since $\MonOlax$ is a $1$-full $2$-subcategory of $\CocarttOlax$, all we need to show is that the right adjoint is in $\MonOlax$, whence the unit and counit are also automatically there.
	Namely, we need to know that $R\colon \DD \to \CC$ is lax $\OO$-monoidal, which follows from the construction, as explained in \cite[Corollary 7.3.2.7]{HA}.
\end{proof}

\subsection{Tensored Categories}

As explained in \cite[Lemma 3.54]{Hei}, the category $\wLMod^\lax_\VV$ is enhanced to a $2$-category, and thus so is its full subcategory $\LMod^\lax_\VV$.
We repeat the construction in more detail, building on the definitions above.
See also \cite[Definition 2.14]{Sha} for a closely related discussion.

We recall the operad $\Ass$ classifying an associative algebra, and the operad $\LM$ classifying an associative algebra and a left module over it.
Consider the monoidal category $\VV$ as an object in $\Mon(\Catt) \cong \Mon_{\Ass} \subset \Mon_{\Ass}^\lax$.
This leads us to the following:

\begin{defn}\label{lmod-lax}
	We define the $2$-category of \tdef{$\VV$-tensored categories and lax $\VV$-linear functors} to be $\mdef{\LModtlax_\VV} := \{\VV\} \times_{\Mon_{\Ass}^\lax} \Mon_{\LM}^\lax$.
	Similarly, we denote by $\LMODtlax_\VV$ the (huge) $2$-category of large $\VV$-tensored categories and lax $\VV$-linear functors.
\end{defn}

\begin{defn}
	We define the $2$-category of \tdef{closed $\VV$-tensored categories and lax $\VV$-linear functors} to be the full $2$-subcategory $\mdef{\LModt^{\cl,\lax}_\VV} \subset \LModtlax_\VV$ on those $\VV$-tensored categories on which the $\VV$-action is closed, and similarly for the large version $\LMODt^{\cl,\lax}_\VV$.
\end{defn}

As a particular case of \cref{lax-O-adj}, we deduce the following, which is an ``if and only if" version of \cite[Example 7.3.2.8 and Remark 7.3.2.9]{HA} phrased in an appropriate $2$-category.

\begin{cor}\label{lax-lin-adj}
	The left adjoints in $\LModtlax_\VV$ are those lax $\VV$-linear functors that are strong $\VV$-linear and left adjoint (on the underlying category)
	$$
	(\LModtlax_\VV)^\LL = \LModt^\fwL_\VV.
	$$
\end{cor}

\begin{proof}
	The result follows from \cref{lax-O-adj} by taking $\OO = \LM$, and pulling back over $\VV$.
	In more detail, consider the following diagram:
	\[\begin{tikzcd}
		{(\LModtlax_\VV)^\LL} & {(\Mon_{\LM}^\lax)^\LL} \\
		{\{\VV\}^\LL} & {(\Mon_{\Ass}^\lax)^\LL} \\
		{\{\VV\}} & {\Mon_{\Ass}^\lax}
		\arrow[from=1-1, to=2-1]
		\arrow["\wr"', from=2-1, to=3-1]
		\arrow[from=3-1, to=3-2]
		\arrow[from=2-2, to=3-2]
		\arrow[from=1-2, to=2-2]
		\arrow[from=1-1, to=1-2]
		\arrow[from=2-1, to=2-2]
	\end{tikzcd}\]
	The bottom square is a pullback square because the lower right morphism is an inclusion of a $2$-subcategory, and the lower left morphism is the map from a point to itself.
	The upper square is a pullback square, because $(-)^\LL$ commutes with limits by \cref{L-lims}.
	Therefore, the outer square is also a pullback square.
	We thus finish by \cref{lax-O-adj}.
\end{proof}

\begin{remark}
	\cref{lmod-lax} and \cref{lax-lin-adj} work for an arbitrary monoidal category $\VV$, without needing to assume that it is \emph{presentably} monoidal.
\end{remark}

We now enhance $\PrLV$ of \cref{prlv} to a $2$-category.

\begin{defn}
	We define the $2$-category of \tdef{presentably $\VV$-tensored categories} to be the $1$-full $2$-subcategory
	$\mdef{\PrtLV} := \LModt_\VV(\PrtL) \subset \LMODtlax_\VV$
	on the objects and morphisms of $\PrLV$.
\end{defn}

\begin{prop}\label{prlv-lmod-l}
	There is a full $2$-subcategory inclusion
	$$
	\PrtLV \subset (\LMODt^{\cl,\lax}_\VV)^\LL.
	$$
\end{prop}

\begin{proof}
	First, we claim that there is a full $2$-subcategory inclusion
	$$
	\PrtLV = \LModt_\VV(\PrL) \subset \LMODt^{\cl,\fwL}_\VV.
	$$
	It is indeed an inclusion, since a presentably $\VV$-tensored category is automatically closed, and morphism in both categories are the $\VV$-linear functors that are left adjoint on the underlying category.
	Second, by \cref{lax-lin-adj}, the target is $(\LMODt^{\cl,\lax}_\VV)^\LL$, which finishes the argument.
\end{proof}

\subsection{Tensored Categories to Enriched Categories}

We recall that $\CatV$ can be enhanced into a $2$-category, as explained for example in \cite[Construction 8.2]{Hei-arxiv}.
Indeed, by \cite[Proposition 5.7.16]{GH}, the functor $\Alg(\PrL) \to \PrL$ given by $\VV \mapsto \CatV$ is lax symmetric monoidal.
Since $\VV$ is a $\Spaces$-module in $\PrL$, we get that $\CatV$ is a $\Cat$-module in $\PrL$ as well, which indeed models a $2$-category (for example via $\rchi$).

\begin{defn}
	We denote the $2$-category of $\VV$-enriched categories by $\mdef{\CattV}$.
\end{defn}

In \cite[Theorem 8.3]{Hei-arxiv}, Heine shows that $\rchi$ of \cref{chi-w} enhances to a $2$-equivalence
$$
\rchi\colon \wLModt^{\cl,\lax}_\VV \iso \CattV.
$$

Using this we deduce the following:

\begin{cor}\label{chi-2}
	The functor from \cref{chi} enhances to a $2$-fully faithful $2$-functor
	$$
	\rchi\colon \PrtLV \to (\CATtV)^\LL.
	$$
\end{cor}

\begin{proof}
	Consider the large version of the $2$-equivalence above, namely
	$\wLMODt^{\cl,\lax}_\VV \iso \CATtV$,
	and restrict it to the full $2$-subcategory
	$\LMODt^{\cl,\lax}_\VV$.
	Taking left adjoints and using \cref{prlv-lmod-l}, we get the $2$-fully faithful $2$-functor
	$
	\PrtLV
	\subset (\LMODt^{\cl,\lax}_\VV)^\LL
	\to (\CATtV)^\LL
	$.
\end{proof}

\subsection{Evaluation and Enriched Hom in Tensored Categories}

In this subsection, we introduce $\hom^\VV$ in the presentably $\VV$-tensored context and show that it satisfies expected properties for adjunction in \cref{lin-adj} (see also \cref{chi-enrch} for the connection to adjunction in $\VV$-enriched categories).
With this definition in mind, we finish the proof of \ref{yon-lem} of \cref{hin-yon} in \cref{yon-lem-restate}.

Let $\CC \in \PrLV$.
We have an equivalence of categories $\Fun^\LL_\VV(\VV, \CC) \iso \CC$ given by evaluation at $\unit_\VV$.
Recall from \cref{prlv-lmod-l} that there is a full $2$-subcategory inclusion $\PrtLV \subset (\LMODt^{\cl,\lax}_\VV)^\LL$.
Thus, passing to right adjoint, gives us a map $\Fun^\LL_\VV(\VV, \CC) \to \Fun_{\laxV}(\CC, \VV)^\op$.
Composing the two and taking $(-)^\op$, we get $\CC^\op \to \Fun_{\laxV}(\CC, \VV)$.

\begin{defn}\label{hom-v-def}
	We let
	$\mdef{\hom^\VV(-, -)}\colon \CC^\op \times \CC \to \VV$
	to be the functor corresponding to the above under the exponential adjunction.
	By construction, it is lax $\VV$-linear in the second argument.
\end{defn}

Namely, this construction sends $X$ to the lax $\VV$-linear functor $\hom^\VV(X, -)\colon \CC \to \VV$, right adjoint to the $\VV$-linear left adjoint $- \otimes X\colon \VV \to \CC$.

\begin{prop}\label{lin-adj}
	Let $L\colon \CC \to \DD \in \PrLV$ be a $\VV$-linear left adjoint functor, and let $R\colon \DD \to \CC$ denote its right adjoint.
	Then there is a natural isomorphism
	$$
	\hom^\VV(L(-), -) \cong \hom^\VV(-, R(-))
	$$
	of functors $\CC^\op \times \DD \to \VV$ lax $\VV$-linear in the second coordinate.
\end{prop}

\begin{proof}
	Indeed, consider the following diagram:
	\[\begin{tikzcd}[cramped]
		\CC & {\Fun^\LL_\VV(\VV, \CC)} & {\Fun_{\laxV}(\CC, \VV)^\op} \\
		\DD & {\Fun^\LL_\VV(\VV, \DD)} & {\Fun_{\laxV}(\DD, \VV)^\op}
		\arrow["L", from=1-1, to=2-1]
		\arrow["{L \circ-}", from=1-2, to=2-2]
		\arrow["{- \circ R}", from=1-3, to=2-3]
		\arrow["\sim", from=2-1, to=2-2]
		\arrow[from=2-2, to=2-3]
		\arrow["\sim", from=1-1, to=1-2]
		\arrow[from=1-2, to=1-3]
	\end{tikzcd}\]
	The left square commutes by passing to the right adjoints and noting that evaluation at $\unit_\VV$ commutes with post-composition.
	To see that the right square commutes, recall that the horizontal morphisms were constructed by the full $2$-subcategory inclusion $\PrtLV \subset (\LMODt^{\cl,\lax}_\VV)^\LL$ and passing to the right adjoints, and the vertical morphisms are also adjoints in the same category, so the commutativity is the fact that the right adjoint of a composition is the composition of the right adjoints in reverse order.
\end{proof}

\begin{remark}\label{chi-enrch}
	We note that \cref{lin-adj} is in the context of presentably $\VV$-tensored categories, and not that of $\VV$-enriched categories.
	One of the main features of $\rchi$ is that the $\VV$-enrichment of $\rchi(\CC)$ is given by $\hom^\VV(X, Y)$, as is shown in \cite[Corollary 7.4.9]{GH}.
	By \cref{chi-2}, the adjunction $L \dashv R$ produces an adjunction $\rchi(L) \dashv \rchi(R)$ in $\CATtV$.
	One might expect that this adjunction provides a natural isomorphism similar to the one of \cref{lin-adj} for the $\VV$-enriched $\hom$'s of $\rchi(\CC)$ and $\rchi(\DD)$.
	However, we do not provide such a natural isomorphism, nor show its compatibility with the one constructed above.
\end{remark}

Let $\CC_0 \in \CatV$, $X \in \CC_0$ and $\DD \in \PrLV$.
Recall from \cref{fun-model} that $\Fun_\VV(\CC_0, \DD) := \LMod_{\CC_0}(\Fun(\CC_0^\simeq, \DD))$.
Consider $i\colon \pt \to \CC_0^\simeq$ choosing $X$, and the following diagram:
\[\begin{tikzcd}
	\DD & {\Fun(\pt, \DD)} && {\Fun(\CC_0^\simeq, \DD)} && {\LMod_{\CC_0}(\Fun(\CC_0^\simeq, \DD))}
	\arrow["{i_!}", shift left=1, dashed, from=1-2, to=1-4]
	\arrow["{i^*}", shift left=1, from=1-4, to=1-2]
	\arrow["\free", shift left=1, dashed, from=1-4, to=1-6]
	\arrow["\forget", shift left=1, from=1-6, to=1-4]
	\arrow[Rightarrow, no head, from=1-2, to=1-1]
\end{tikzcd}\]
Here the dashed arrows are the left adjoints of the solid arrows.
Also recall from \cref{pshv-def} that if we replace $\CC_0$ by $\CC_0^{\op}$ and let $\DD = \VV$, the $\VV$-functor category is $\PShV(\CC_0)$.
Moreover, in this case the three categories in the diagram are presentably $\VV$-tensored and the two left adjoint functors are $\VV$-linear, as in \cite[6.2.3]{Hin}.
Via \cref{prlv-lmod-l}, the right adjoints are canonically lax $\VV$-linear.

\begin{defn}\label{eval-x}
	Let $\CC_0 \in \CatV$, $X \in \CC_0$ and $\DD \in \PrLV$ be as above.
	We define the \tdef{evaluation at $X$} functor
	$\mdef{\eval_X}\colon \Fun_\VV(\CC_0, \DD) \to \DD$
	by $\eval_X := i^* \circ \forget$.
	For the case $\eval_X\colon \PShV(\CC_0) \to \VV$, it is canonically lax $\VV$-linear.
\end{defn}

\begin{prop}\label{colim-levelwise}
	Let $\CC_0 \in \CatV$ and $\DD \in \PrLV$, and consider $\Fun_\VV(\CC_0, \DD) \in \PrL$.
	Then $\eval_X$ commutes with (co)limits (i.e.\ (co)limits are computed level-wise), and are jointly conservative over all $X$.
	For $\PShV(\CC_0) \in \PrLV$, the lax $\VV$-linear structure on $\eval_X$ is strong (i.e., the $\VV$-action is level-wise).
\end{prop}

\begin{proof}
	Recall that $\eval_X = i^* \circ \forget$.
	The functor $\forget$ commutes with all (co)limits by \cite[Corollary 4.2.3.3 and Corollary 4.2.3.5]{HA} as the forgetful from modules, while $i^*$ commutes with (co)limits as they are computed level-wise in (unenriched) functor categories.

	The forgetful from modules is always conservative by \cite[Corollary 4.2.3.2]{HA}, and the evaluation at $X$ functors $\Fun(\CC_0^\simeq, \DD) \to \DD$ are jointly conservative.

	For the case of $\PShV(\CC_0)$, the $\VV$-action is by construction given level-wise, as explained in \cite[6.2.3]{Hin}.
\end{proof}

Finally, we are in position to prove our variant of Hinich's enriched Yoneda lemma, appearing as \ref{yon-lem} of \cref{hin-yon}.

\begin{prop}\label{yon-lem-restate}
	Let $\CC_0 \in \CatV$, then $\hom^\VV(\yonV(X), -) \cong \eval_X$ as lax $\VV$-linear functors, and in particular $\hom^\VV(\yonV(X), -)$ is also strong $\VV$-linear.
\end{prop}

\begin{proof}
	We need to show that $\hom^\VV(\yonV(X), -) \cong \eval_X$ as lax $\VV$-linear functors.
	By construction, this is equivalent to showing that $- \otimes \yonV(X) \cong \free \circ i_!$ as $\VV$-linear functors.
	Consider the following diagram:
	\[\begin{tikzcd}
		\VV &&& {\Fun(\CC_0^{\op,\simeq}, \VV)} && {\LMod_{\CC_0^\op}(\Fun(\CC_0^{\op,\simeq}, \VV))} \\
		\\
		\Spaces &&& {\Fun(\CC_0^{\op,\simeq}, \Spaces)}
		\arrow["{-\otimes \unit_\VV}", shift left=1, dashed, from=3-1, to=1-1]
		\arrow["{\hom(\unit_\VV,-)}", shift left=1, from=1-1, to=3-1]
		\arrow["{-\otimes \unit_\VV}", shift left=1, dashed, from=3-4, to=1-4]
		\arrow["{\hom(\unit_\VV,-)}", shift left=1, from=1-4, to=3-4]
		\arrow["{i_!}", shift left=1, dashed, from=3-1, to=3-4]
		\arrow["{i^*}", shift left=1, from=3-4, to=3-1]
		\arrow["{i_!}", shift left=1, dashed, from=1-1, to=1-4]
		\arrow["{i^*}", shift left=1, from=1-4, to=1-1]
		\arrow["\free", shift left=1, dashed, from=1-4, to=1-6]
		\arrow["\forget", shift left=1, from=1-6, to=1-4]
	\end{tikzcd}\]
	Here the dashed arrows are the left adjoints of the solid arrows.
	Clearly, the solid square commutes, and thus, the dashed square obtained by passing to left adjoints also commutes.
	
	By \cite[6.2.6]{Hin}, we have that $\yonV(X) \cong \free((- \otimes \unit_\VV) \circ \yon(X))$, where here $\yon(X) \in \Fun(\CC_0^{\op,\simeq}, \Spaces)$ is the image of $X$ under the (unenriched) Yoneda embedding.
	Applying the naturality of the (unenriched) Yoneda embedding to $i\colon \pt \to \Spaces$, and noting that $\yon\colon \pt \to \Spaces$ sends $\yon(\pt) = \pt$, we get that $\yon(X) = \yon(i(\pt)) \cong i_!(\yon(\pt)) \cong i_!(\pt)$.
	Using the commutativity of the dashed square above, we conclude that
	$$
	\yonV(X)
	\cong \free((- \otimes \unit_\VV) \circ \yon(X))
	\cong \free((- \otimes \unit_\VV) \circ i_!(\pt))
	\cong \free(i_! \unit_\VV)
	.
	$$
	Since $\free$ and $i_!$ are $\VV$-linear functors, we finally get that
	$$
	- \otimes \yonV(X)
	\cong - \otimes \free(i_! \unit_\VV)
	\cong \free \circ i_!(- \otimes \unit_\VV)
	\cong \free \circ i_!
	$$ 
	as $\VV$-linear functors, concluding the proof.
\end{proof}

%% file: atomics.tex
\section{Atomic Objects, Presheaves and Yoneda}\label{sec-at}

In this section, building on the $2$-categorical results of the previous section, we study internally left adjoint functors and atomic objects, and connect them to the enriched Yoneda embedding.

We begin in \cref{iL} by defining internally left adjoints between presentably $\VV$-tensored categories.
Namely, $\VV$-linear left adjoint functors, whose right adjoints are also left adjoint and their canonical lax $\VV$-linear structure is strong.
Following this, in \cref{at-def} we define a finiteness condition called being atomic.
We say that $X$ in $\CC \in \PrLV$ is atomic if $- \otimes X\colon \VV \to \CC$ is internally left adjoint, that is, if $\hom^\VV(X, -)$ commutes with colimits and is $\VV$-linear.
For instance, in the case $\VV = \Sp$, this coincides with condition of being compact.
From the definition, internally left adjoints send atomic objects to atomic objects.
The key technical result of this section is \cref{atomic-pres-il}, giving a converse result under the assumption that the source category is generated from the atomics under weighted colimits.

Next, in \cref{yon-at} we show that $\yonV(X)$ is atomic in $\PShV(\CC_0)$, and in \cref{psh-molecular} we show that together they generate $\PShV(\CC_0)$ under weighted colimits.
This allows us to use \cref{atomic-pres-il} to prove \cref{yon-at-adj}, the main result of this section.
This result says that the partial adjunction of \cref{pshv-chi-adj} restricts to a (non-partial) adjunction between the enriched presheaves functor and taking the atomics, with the unit being (the factorization through the atomics of) the enriched Yoneda embedding.

\subsection{Internally Left Adjoints}

\begin{defn}\label{iL}
	A $\VV$-linear left adjoint functor is called \tdef{internally left adjoint} if it is left adjoint in the $2$-category $\PrtLV$.
	We denote the wide $2$-subcategory on the internally left adjoint functors by $\mdef{\PrtiLV} := (\PrtLV)^\LL$.
	For $\CC, \DD \in \PrtiLV$, we denote the category of $\VV$-linear internally left adjoint functors between them by $\mdef{\Fun^\intL_\VV(\CC, \DD)} \in \CAT$ and the corresponding $\mdef{\hom^\intL_\VV(\CC, \DD)} \in \SPACES$.
\end{defn}

Note that for any $\VV$-linear left adjoint functor $L\colon \CC \to \DD$, \cref{prlv-lmod-l} shows that it admits a lax $\VV$-linear right adjoint $R\colon \DD \to \CC$.

\begin{prop}\label{iL-cond}
	A $\VV$-linear left adjoint functor $L\colon \CC \to \DD$ is internally left adjoint if and only if the lax $\VV$-linear right adjoint $R\colon \DD \to \CC$ is strong and itself a left adjoint.
\end{prop}

\begin{proof}
	As above, $R$ is a morphism in $\LMODt^{\cl,\lax}_\VV$ between objects of the $2$-subcategory $\PrtLV$, and the condition is that it is in fact a morphism in that $2$-subcategory.
	By \cref{prlv-lmod-l}, $\PrtLV$ is a full $2$-subcategory of $(\LMODt^{\cl,\lax}_\VV)^\LL$, thus the condition is that $R$ is a morphism in that $2$-subcategory.
	The result then follows from \cref{lax-lin-adj}.
\end{proof}

Taking left adjoints in \cref{chi-2}, we immediately get:

\begin{cor}\label{chi-2-l}
	$\rchi$ restricts to a $2$-fully faithful $2$-functor
	$$
	\rchi\colon \PrtiLV \to (\CATtV)^{\LL\LL}.
	$$
\end{cor}

\subsection{Atomic Objects}

\begin{defn}\label{at-def}
	Let $\CC \in \PrLV$.
	We say that $X \in \CC$ is \tdef{atomic} if the functor
	$- \otimes X\colon \VV \to \CC$
	is internally left adjoint.
	We denote the full $\VV$-subcategory on the atomic objects by $\mdef{\CC^{\atomic}} \subset \rchi(\CC)$.
\end{defn}

\begin{example}
	The unit $\unit_\VV \in \VV$ is always atomic, because $- \otimes \unit_\VV$ is the identity functor.
\end{example}

\begin{example}[{\cite[Proposition 2.8]{BMS}}]
	In the case $\VV = \Sp$, atomic objects coincide with compact objects.
\end{example}

Recall that for $\CC \in \PrLV$, there is an equivalence
$\Fun^\LL_\VV(\VV, \CC) \iso \CC$
given by evaluation at $\unit_\VV$, with inverse sending $X$ to $- \otimes X\colon \VV \to \CC$.
Thus, immediately from the definition, we deduce the following:

\begin{prop}
	Evaluation at $\unit_\VV$ induces an equivalence
	$\Fun^\intL_\VV(\VV, \CC) \iso \CC^{\atomic}$.
\end{prop}

Since the right adjoint of $- \otimes X\colon \VV \to \CC$ is $\hom^\VV(X, -)\colon \CC \to \VV$, the following follows immediately from \cref{iL-cond}.

\begin{prop}\label{at-w-colims}
	An object $X \in \CC$ is atomic if and only if $\hom^\VV(X, -)\colon \CC \to \VV$ preserves colimits and the lax $\VV$-linear structure is strong.
\end{prop}

\begin{prop}
	The atomics are a small $\VV$-enriched category.
\end{prop}

\begin{proof}
	The argument is identical to \cite[Proposition 2.9]{BMS}.
	We repeat the details for the convenience of the reader.

	Let $\kappa$ be a regular cardinal such that the unit $\unit_\VV$ is $\kappa$-compact, namely $\hom(\unit_\VV, -)\colon \VV \to \Spaces$ commutes with $\kappa$-filtered colimits.
	We show that all atomic objects are $\kappa$-compact.
	Let $\CC \in \PrLV$ and let $X \in \CC$ be atomic.
	By \cref{at-w-colims}, $\hom^\VV(X, -)\colon \CC \to \VV$ commutes with $\kappa$-filtered colimits.
	Since
	$\hom(X, -) \cong \hom(\unit_\VV \otimes X, -) \cong \hom(\unit_\VV, \hom^\VV(X, -))$,
	we get that $\hom(X, -)$ also commutes with $\kappa$-filtered colimits, i.e.\ $X$ is $\kappa$-compact.
	We have shown that $\CC^{\atomic} \subseteq \CC^\kappa$, the latter being a small category, concluding the proof.
\end{proof}

\begin{prop}\label{il-atomic-pres}
	Let $L\colon \CC \to \DD$ be an internally left adjoint, i.e.\ a morphism in $\PriLV$, then it sends atomics to atomics.
	Thus, the $\VV$-functor $\rchi(L)\colon \rchi(\CC) \to \rchi(\DD)$ factors to a $\VV$-functor $L\colon \CC^{\atomic} \to \DD^{\atomic}$.
\end{prop}

\begin{proof}
	Let $X \in \CC^{\atomic}$.
	Then, since $L$ is $\VV$-linear, we have a natural isomorphism $- \otimes LX \cong L(- \otimes X)$ of functors $\VV \to \DD$.
	Since both $- \otimes X$ and $L$ are internally left adjoints, so is $- \otimes LX$, i.e.\ $LX$ is atomic as required.
\end{proof}

Recall that $\CC^{\atomic}$ is a full $\VV$-subcategory of $\rchi(\CC)$, thus $\CC^{\atomic} \to \rchi(\CC)$ is a subobject (indeed, the space of maps to $\CC^{\atomic}$ is exactly the subspace of maps to $\rchi(\CC)$ that land in $\CC^{\atomic}$).
Furthermore, by \cref{il-atomic-pres}, the restriction of internally left adjoint functors to the atomics factor through the atomics.
Thus, by \cite[Proposition A.1]{Ram}, this assembles into an induced subfunctor, which furthermore lands in $\CatV \subset \CATV$.

\begin{defn}\label{atomic-def}
	We denote the induced subfunctor of \tdef{atomics} by
	$$
	\mdef{(-)^{\atomic}}\colon \PriLV \to \CatV,
	$$
	equipped with a natural transformation
	$(-)^{\atomic} \Rightarrow \rchi|_{\PriLV}$
	of functors $\PriLV \to \CATV$.
\end{defn}

\begin{defn}
	Let $\CC \in \PrLV$.
	We say that a collection of atomic objects $B \subseteq \CC^{\atomic}$ are \tdef{atomic generators}, if $\CC$ is generated from $B$ under weighted colimits.
	If such $B$ exists, we say that $\CC$ is \tdef{molecular}.
\end{defn}

\begin{remark}
	We note that this definition a priori differs from our definition of molecular in \cite{BMS} (where we work over a mode). 
	We expect that closure under weighted colimits coincides with closure under colimits and the $\VV$-action, but we are unaware of a proof of this statement.
	Such a result would make the connection transparent.
\end{remark}

We now wish to prove a converse to \cref{il-atomic-pres} under the assumption that the source is molecular.
To that end, we begin with the following.

\begin{defn}\label{hom-f}
	Let $\II \in \CatV$ and $\CC \in \PrLV$, and fix $f \in \Fun_\VV(\II, \CC)$.
	By \ref{wcol-lin} of \cref{hin-yon}, the functor
	$\coliminl_{\II}^{(-)}(f)\colon \PShV(\II) \to \CC$
	is colimit preserving, i.e.\ a left adjoint, and $\VV$-linear.
	We denote the lax $\VV$-linear right adjoint by
	$$
	\mdef{\hom^\VV(f(-), -)}\colon \CC \to \PShV(\II).
	$$
\end{defn}

\begin{lem}\label{wcol-adj}
	There is a natural isomorphism of lax $\VV$-linear functors
	$$
	\eval_i \circ \hom^\VV(f(-), -) \cong \hom^\VV(f(i), -).
	$$
\end{lem}

\begin{proof}
	Observe that we have an equivalence of $\VV$-linear left adjoint functors
	$$
		\coliminl_{\II}^{(-)}(f) \circ (- \otimes \yonV(i))
		= \coliminl_{\II}^{(- \otimes \yonV(i))}(f)
		\overset{\text{(1)}}{\cong} - \otimes \coliminl_{\II}^{\yonV(i)}(f)
		\overset{\text{(2)}}{\cong} - \otimes f(i),
	$$
	where (1) follows from \ref{wcol-lin} of \cref{hin-yon}, and (2) follows from \ref{psh-univ} of \cref{hin-yon}.
	Passing to the lax $\VV$-linear right adjoints, we get
	$$
	\hom^\VV(\yonV(i), \hom^\VV(f(-), -)) \cong \hom^\VV(f(i), -).
	$$
	We finish by recalling that $\hom^\VV(\yonV(i), -)$ is the evaluation at $i$ by \ref{yon-lem} of \cref{hin-yon}.
\end{proof}

\begin{prop}\label{atomic-pres-il}
	Let $L\colon \CC \to \DD$ be in $\PrLV$.
	If $\CC$ is molecular and $L$ sends a collection of atomics generators $B \subset \CC^{\atomic}$ to atomic objects of $\DD$, then $L$ is internally left adjoint.
\end{prop}

\begin{proof}
	We adapt the proof of \cite[Proposition 2.14]{BMS}.
	We wish to show that $R\colon \DD \to \CC$, the right adjoint of $L$, is itself a left adjoint and that the lax $\VV$-action is strong.
	Thus, it suffices to show that for any $v \in \VV$ and $Y\colon \II \to \DD$, the canonical map $v\otimes\colim_I R Y_i \to R(v\otimes\colim_I Y_i)$ is an isomorphism.
	By the (unenriched) Yoneda lemma in the category $\CC$, this is equivalent to checking that for every $X \in \CC$ the map
	\begin{equation}
		\hom(X, v\otimes\colim_I R Y_i) \to \hom(X, R(v\otimes\colim_I Y_i))
	\end{equation}
	is an isomorphism.
	We will in fact show the stronger statement that
	\begin{equation}
		\hom^\VV(X, v\otimes\colim_I R Y_i) \to \hom^\VV(X, R(v\otimes\colim_I Y_i))
	\end{equation}
	is an isomorphism, which implies the previous statement by taking $\hom(\unit_\VV, -)$.
	Let $A \subset \CC$ be the collection of objects $X$ for which it is an isomorphism, and we will show that $A = \CC$.
	
	We first show that $B \subset A$.
	Let $X \in B$.
	By \cref{lin-adj}, the following diagram commutes, and both vertical maps are isomorphisms:
	$$\begin{tikzcd}
		v\otimes\colim_I \hom^\VV(X, R Y_i) \arrow{r}{} \arrow{d}{\wr}
			& \hom^\VV(X, v\otimes\colim_I R Y_i) \arrow{r}{}
			& \hom^\VV(X, R(v\otimes\colim_I Y_i)) \arrow{d}{\wr}
		\\
		v\otimes\colim_I \hom^\VV(L X, Y_i) \arrow{rr}{}
			&
			& \hom^\VV(L X, v\otimes\colim_I Y_i)
	\end{tikzcd}$$
	The upper-left morphism is an isomorphism because $X$ is atomic, and similarly the bottom morphism is an isomorphism because $LX$ is atomic since $X \in B$ and $L$ sends $B$ to atomic objects.
	This shows that the upper-right morphism is an isomorphism as well.

	We now show that $A$ is closed under weighted colimits.
	Let $J \in \CatV$, $W \in \PShV(J)$, and $f \in \Fun_\VV(J, \CC)$ a $\VV$-functor landing in $A \subset \CC$.
	We show that $X := \coliminl_J^W(f)$ is in $A$.
	By definition, $\hom^\VV(f(-), -)\colon \CC \to \PShV(J)$ is the lax $\VV$-linear right adjoint of $\coliminl_J^{(-)}(f)\colon \PShV(J) \to \CC$.
	Thus, by \cref{lin-adj} we have a natural isomorphism of functors $\CC \to \VV$
	$$
	\hom_\CC^\VV(X, -)
	= \hom_\CC^\VV(\coliminl_J^W(f), -)
	\cong \hom_{\PShV(J)}^\VV(W, \hom^\VV(f(-), -)).
	$$
	Thus it suffices to check that the map
	$$
	\hom_\CC^\VV(f(-), v\otimes\colim_I R Y_i)
	\to \hom_\CC^\VV(f(-), R(v\otimes\colim_I Y_i))
	\quad \in \PShV(J)
	$$
	is an isomorphism.
	By \cref{colim-levelwise}, the evaluation at $j \in J$ are jointly conservative, so it suffices to check that the map is an isomorphism after evaluation at every $j$.
	By \cref{wcol-adj}, this means that we need to check that
	$$
	\hom^\VV(f(j), v\otimes\colim_I R Y_i)
	\to \hom^\VV(f(j), R(v\otimes\colim_I Y_i))
	$$
	is an isomorphism, which holds since by assumption $f$ lands in $A$.

	Recall that $B$ are atomic generators of $\CC$, and we have shown that $B \subset A$ and that $A$ is closed under weighted colimits, thus $A = \CC$, as required.
\end{proof}

\subsection{Atomics--Presheaves Adjunction}

\begin{prop}\label{yon-at}
	Let $\CC_0 \in \CatV$, then the enriched Yoneda embedding lands in the atomic objects, yielding a $\VV$-functor
	$\yonV\colon \CC_0 \to \PShV(\CC_0)^{\atomic}$.
\end{prop}

\begin{proof}
	Recall from \cref{at-w-colims} that $\yonV(X)$ is atomic if and only if $\hom^\VV(\yonV(X), -)\colon \PShV(\CC_0) \to \VV$ preserves colimits and the lax $\VV$-linear structure is strong.
	By \ref{yon-lem} of \cref{hin-yon}, this $\VV$-functor is the evaluation at $X$, concluding the proof by \cref{colim-levelwise}.
\end{proof}

We thus get an induced natural transformation
$\yonV\colon \id \Rightarrow \PShV(-)^{\atomic}$
of functors $\CatV \to \CatV$.

\begin{prop}\label{psh-molecular}
	The category $\PShV(\CC_0)$ is molecular, with the image of the enriched Yoneda embedding as atomic generators.
\end{prop}

\begin{proof}
	By \cref{yon-at}, the image of $\yonV$ are atomic in $\PShV(\CC_0)$.
	The fact that $\PShV(\CC_0)$ is generated from the image of enriched Yoneda embedding under weighted colimits is \cite[6.3.1]{Hin2}, recalled as \ref{psh-univ} of \cref{hin-yon}.
\end{proof}

\begin{prop}\label{pshv-il}
	The functor $\PShV\colon \CatV \to \PrLV$ of \cref{pshv-chi-adj} lands in $\PriLV$.
\end{prop}

\begin{proof}
	Let $f\colon \CC_0 \to \DD_0$ be a $\VV$-functor, and we need to show that $f_!\colon \PShV(\CC_0) \to \PShV(\DD_0)$ is internally left adjoint.
	Recall from \cref{psh-molecular} that $\PShV(\CC_0)$ is molecular with the image of $\yonV_{\CC_0}$ as atomic generators.
	As in \cref{yon-nat}, naturality of the enriched Yoneda embedding says that $\yonV_{\DD_0} f \cong f_! \yonV_{\CC_0}$.
	Thus, $f_!$ sends the image of $\yonV_{\CC_0}$ to the image of $\yonV_{\DD_0}$ which are atomic in $\PShV(\DD_0)$.
	Thus, $f_!$ sends a collection of atomic generators to atomic objects, so it is internally left adjoint by \cref{atomic-pres-il}.
\end{proof}

\begin{thm}\label{yon-at-adj}
	The partial adjunctions of \cref{pshv-chi-adj} restricts to an adjunction
	$$
	\PShV\colon \CatV \rightleftarrows \PriLV\noloc (-)^{\atomic}
	$$
	with unit the enriched Yoneda embedding $\yonV\colon \id \Rightarrow \PShV(-)^{\atomic}$.
\end{thm}

\begin{proof}
	We need to check that for any $\CC_0 \in \CatV$ and $\DD \in \PriLV$, the map
	\begin{equation}\label{pv-at-check}
		\hom^\intL_\VV(\PShV(\CC_0), \DD)
		\xrightarrow{(-)^{\atomic}} \hom_\VV(\PShV(\CC_0)^{\atomic}, \DD^{\atomic})
		\xrightarrow{(\yonV)^*} \hom_\VV(\CC_0, \DD^{\atomic})
	\end{equation}
	is an equivalence.
	Recall that
	\begin{equation}\label{pv-at-lift}
		\hom^\LL_\VV(\PShV(\CC_0), \DD)
		\xrightarrow{\rchi} \hom_\VV(\rchi(\PShV(\CC_0)), \rchi(\DD))
		\xrightarrow{(\yonV)^*} \hom_\VV(\CC_0, \rchi(\DD))
	\end{equation}
	is an equivalence.
	Furthermore, both the first and last spaces in \eqref{pv-at-check} are a collection of connected components of the first and last spaces in \eqref{pv-at-lift}, showing that the composition in \eqref{pv-at-check} is an inclusion of connected components.

	To finish the argument, we need to show that the composition in \eqref{pv-at-check} hits every connected component.
	To that end, let $f\colon \CC_0 \to \DD^{\atomic}$ be a $\VV$-functor.
	We can post-compose it with the inclusion $\VV$-functor $\DD^{\atomic} \to \DD$, and using \eqref{pv-at-lift} we get $\tilde{f}\colon \PShV(\CC_0) \to \DD$ in $\PrLV$.
	It is left to show that $\tilde{f}$ is internally left adjoint.
	Recall that the image of the enriched Yoneda embedding forms a collection of atomic generators by \cref{psh-molecular}.
	Furthermore, by construction, $\tilde{f}(\yonV(X)) = f(X) \in \DD^{\atomic}$ is atomic.
	Thus, we have shown that $\tilde{f}$ sends a collection of atomic generators to atomics, so it is indeed internally left adjoint by \cref{atomic-pres-il}.
\end{proof}

We extend our notations from \cref{pshv-of-f}:

\begin{defn}\label{pshv-of-f-2}
	For $f\colon \CC_0 \to \DD_0$ a morphism in $\CatV$, consider the internally left adjoint functor $f_!\colon \PShV(\CC_0) \to \PShV(\DD_0)$.
	It has a right adjoint in $\PrtLV$, i.e.\ a $\VV$-linear left and right adjoint functor denoted $\mdef{f^\circledast}\colon \PShV(\DD_0) \to \PShV(\CC_0)$.
	This functor thus has a further lax $\VV$-linear right adjoint denoted $\mdef{f_\circledast}\colon \PShV(\CC_0) \to \PShV(\DD_0)$.
\end{defn}

Using \cref{chi-2-l} we get the composition
$$
\rchi\PShV\colon \CatV \xrightarrow{\PShV} \PriLV \xrightarrow{\rchi} (\CATV)^{\LL\LL},
$$
allowing to pass the (double) adjunction to $\CATtV$.

\begin{cor}
	For $f\colon \CC_0 \to \DD_0$ a morphism in $\CatV$, we get a double adjunction $\rchi(f_!) \dashv \rchi(f^\circledast) \dashv \rchi(f_\circledast)$ in $\CATtV$.
\end{cor}

\begin{cor}\label{upperstar-pointwise}
	Let $f\colon \CC_0 \to \DD_0$ in $\CatV$, then $f^\circledast(G)(X) \cong G(f(X))$.
\end{cor}

\begin{proof}
	By \cref{yon-nat} we have $f_! \yonV \cong \yonV f$ in $\Fun_\VV(\CC_0, \PShV(\DD_0))$.
	We thus get
	$$
	f^\circledast(G)(X)
	\cong \hom^\VV(\yonV(X), f^\circledast(G))
	\cong \hom^\VV(f_!(\yonV(X)), G)
	\cong \hom^\VV(\yonV(f(X)), G)
	\cong G(f(X)),
	$$
	where the first and last steps follow from the enriched Yoneda lemma of \ref{yon-lem} of \cref{hin-yon}, the second step is by \cref{lin-adj}, and the third step is by the isomorphism above.
\end{proof}

%% file: hei.tex
\section{Heine's Enriched Yoneda Embedding}\label{sec-hei}

Independently of Hinich's enriched Yoneda embedding, Heine defines an enriched Yoneda embedding as well, which satisfies the exact same universal property appearing in \cref{psh-chi-pre-adj}.

\begin{thm}[{\cite[Theorem 1.11, Definition 5.3]{Hei}}]\label{hei-yon}
	Let $\CC_0 \in \CatV$.
	There is a $\VV$-natural transformation $\HeiyonV\colon \CC_0 \to \rchi(\PShV(\CC_0))$.
	For every $\DD \in \PrLV$, the composition
	$$
	\hom^\LL_\VV(\PShV(\CC_0), \DD)
	\xrightarrow{\rchi} \hom_\VV(\rchi(\PShV(\CC_0)), \rchi(\DD))
	\xrightarrow{(\HeiyonV)^*} \hom_\VV(\CC_0, \rchi(\DD))
	$$
	is an equivalence.
\end{thm}

We remark that this does \emph{not} show that Hinich's and Heine's enriched Yoneda embeddings agree as functors.
Rather, they are identified only up to an automorphism, that is, there exist some automorphism $\psi_{\CC_0}\colon \PShV(\CC_0) \iso \PShV(\CC_0)$ and $\psi_{\CC_0} \yonV \cong \HeiyonV$.
The question of whether $\psi_{\CC_0}$ is the identity is equivalent to the question of whether the two versions of the enriched Yoneda embedding coincide (see also \cref{hei-hin}).

We now explain that the main results of this paper hold for Heine's enriched Yoneda embedding as well, which in particular makes $\psi_{\CC_0}$ natural in $\CC_0$.
We begin by observing that the two versions of the enriched Yoneda embedding agree point-wise.
In particular, Heine's version also satisfies the enriched Yoneda lemma.

\begin{prop}\label{yons-pointwise}
	For every $X \in \CC_0$, there is an isomorphism $\HeiyonV(X) \cong \yonV(X)$.
	In particular, $\hom^\VV(\HeiyonV(X), -)\colon \PShV(\CC_0) \to \VV$ is given by evaluation at $X$.
\end{prop}

\begin{proof}
	Consider the composition
	$$
	\CC_0^\simeq
	\xrightarrow{\yon} \Fun(\CC_0^{\op, \simeq}, \Spaces)
	\xrightarrow{\unit_\VV \otimes -} \Fun(\CC_0^{\op, \simeq}, \VV)
	\xrightarrow{\free} \LMod_{\CC_0^\op}(\Fun(\CC_0^{\op, \simeq}, \VV))
	= \PShV(\CC_0).
	$$
	By \cite[6.2.6]{Hin} and \cite[Discussion under Remark 5.6]{Hei}, $\yonV(X)$ and $\HeiyonV(X)$ (respectively) are the image of $X$ under this composition.

	The second part then follows from \ref{yon-lem} of \cref{hin-yon}.
\end{proof}

Note that this shows that for every $X \in \CC_0$, the automorphism $\psi_{\CC_0}$ is the identity on $\yonV(X)$, but this does not imply that $\psi_{\CC_0}$ is the identity functor on the image of $\yonV$.

We move on to extending the main results of the paper to Heine's enriched Yoneda embedding.
The proofs in \cref{sec-part-adj} and \cref{sec-at} have relied on the universal property of \cref{psh-chi-pre-adj} and the enriched Yoneda lemma of \ref{yon-lem} of \cref{hin-yon} for $\yonV$.
These two results also hold for $\HeiyonV$ by the above.
We also used weighted colimits, and their relationship with the enriched Yoneda embedding appearing in \ref{wcol-lin} and \ref{psh-univ} of \cref{hin-yon}.
However, this was only used in the proof of \cref{psh-molecular}, whose statement is about the image of the enriched Yoneda embedding which is the same for $\yonV$ and $\HeiyonV$ by \cref{yons-pointwise}, and in the proof of \cref{wcol-adj}, whose statement does not involve the enriched Yoneda embedding.
Thus, we see that the results of \cref{sec-part-adj} and \cref{sec-at} hold for Heine's enriched Yoneda embedding as well.
Notably, the analogue of \cref{yon-at-adj} gives an adjunction
$$
\HeiPShV\colon \CatV \rightleftarrows \PriLV\noloc (-)^{\atomic}
$$
with unit Heine's enriched Yoneda embedding $\HeiyonV\colon \id \Rightarrow \HeiPShV(-)^{\atomic}$.

The uniqueness of left adjoints shows that there is a natural isomorphism $\psi\colon \PShV \to \HeiPShV$ compatible with the unit maps, i.e., with the two versions of the enriched Yoneda embedding.
Since by construction $\HeiPShV$ is point-wise given by $\PShV(\CC_0)$, evaluating $\psi$ at $\CC_0$ reproduces $\psi_{\CC_0}\colon \PShV(\CC_0) \to \HeiPShV(\CC_0) = \PShV(\CC_0)$.